\documentclass[10pt]{amsart}
\usepackage{amsmath,amsfonts,amssymb,amsxtra,latexsym,amscd,enumerate,amsthm}

\usepackage{latexsym}
\usepackage{epsfig}

\def\z{\zeta}
\def\t{\theta}
\def\r{\rho}
\def\g{\gamma}
\def\a{\alpha}
\def\q{\frac{1}{2}}

\def\b{\beta}
\def\e{\varepsilon}

\def\l{\lambda}
\def\D{\Delta}
\def\to{\rightarrow}
\def\Ra{\Rightarrow}
\def\les{\lesssim}

\def\R{\mathbb{R}}
\def\C{\mathbb{C}}
\def\Z{\mathbb{Z}}
\def\N{\mathbb{N}}

\def\beq{\begin{equation}}
\def\eeq{\end{equation}}

\def\beq{\begin{equation}}
\def\eeq{\end{equation}}

\newtheorem{t1}{Theorem}
\newtheorem{l5}{Lemma}

\newtheorem{p3}{Proposition}

\newtheorem{r3}{Remark}

\begin{document}
\title[Quadratic Nonlinear Derivative Schr\"odinger Equations]{Quadratic Nonlinear Derivative Schr\"odinger Equations - Part 2}

\author{Ioan Bejenaru}
\address{Department of Mathematics, UCLA, Los Angeles CA 90095-1555}
\email{bejenaru@math.ucla.edu}

\begin{abstract} 

In this paper we consider the local well-posedness theory for the quadratic nonlinear Schr\"odinger equation with low regularity initial data in the case when the nonlinearity contains derivatives. We work in $2+1$ dimensions and prove a local well-posedness result close to scaling for small initial data.

\end{abstract}

\maketitle

\section{Introduction}

This work is concerned with  the initial value problem for the nonlinear Schr\"odinger equations which generically have the form:

\begin{equation} \label{E:s}
\begin{cases}
\begin{aligned} 
iu_{t}-\Delta u &= P(u,\bar{u},\nabla u, \nabla \bar{u}), \   t\in \mathbb{R}, x\in \mathbb{R}^{n} \\
u(x,0) &=u_{0}(x)
\end{aligned}
\end{cases}
\end{equation}

\noindent
where $u: \mathbb{R}^{n} \times \mathbb{R} \rightarrow \mathbb{C}$ and $P: \mathbb{C}^{2n+2} \rightarrow \mathbb{C}$ is a polynomial. 

We are interested in the theory of local well-posedness for this problem in Sobolev spaces. In Part 1, see \cite{be}, we considered the same problem and it would be useful to read the Introduction there. We summarize it in what follows.

We motivated the fact that the problem becomes more difficult once we consider quadratic and higher order nonlinearities. In this case the most general result known is due to Kenig, Ponce and Vega, see \cite{k2}: 

\begin{t1}

Assume that $P$ has no constant or linear terms. Then there exist $s=s(n,P) > 0$ and $m=m(n,P) >0$ such that $\forall u_{0} \in H^{s}(\R^{n}) \cap L^{2}(\R^{n}:|x|^{2m}dx)$ the problem (\ref{E:s}) has a unique solution in $C([0,T]: H^{s} \cap L^{2}(\R^{n}:|x|^{2m}dx)$ where $T=T(||u||_{H^{s} \cap L^{2}(\R^{n}:|x|^{2m}dx) })$.
\end{t1}

 If $P$ does not contain quadratic terms, then above authors also obtain a similar result without involving any decay, see  \cite{k2}.

We outlined the fact that the case when the nonlinearity contains derivatives is more delicate. One of the reasons is the loss of derivative on the right hand side of the equation. The other one is the need of some decay on the initial data. This is motivated by an early result due to Mizohata, see \cite{m}, which proves that for the problem:

\begin{equation}
\begin{cases}
\begin{aligned} 
iu_{t}-\Delta u &= b_{1}(x) \nabla u , \   t\in \mathbb{R}, x\in \mathbb{R}^{n} \\u(x,0) &=u_{0}(x)
\end{aligned}
\end{cases}
\end{equation}

the following condition on $b_{1}$ is necessary for the $L^{2}$ well-posedness theory:

\beq \label{co1}
\sup_{x \in \R^{n}, \omega \in \mathbb{S}^{n-1}, R > 0}  |Re \int_{0}^{R} b_{1}(x+r \omega) \cdot \omega dr| < \infty
\eeq

We also remarked that the use of decay of type $L^{2}(|x|^{m}dx)$ is not the most appropriate for the Schr\"odinger equation since this structure is not conserved under the linear flow.  

Then we stated the goal of the paper. We wanted to know what is the lowest Sobolev regularity the initial data can have so that we have well-posedness? When asking this question, one should be more specific about the type of the equation and the dimension of the space. 
 
The quadratic terms in $P$ are the first ones to be understood. The quadratic nonlinearities without derivatives have been studied in \cite{c1} and the results obtained are close to scaling.

If the nonlinearity contains terms with derivatives then the problem is called derivative non-linear Schr\"odinger equation (D-NLS). 
The results for quadratic (D-NLS) did not yet reach this level of precision, the main difficulty being generated by the loss of one-derivative in the nonlinearity. The scaling exponents for the problem are $s_{c}=\frac{n}{2}-1$, when only one of the terms contains derivatives (for instance $u Du$), and $s_{c}=\frac{n}{2}$, when both terms contain contain derivatives (for instance $Du Du$). The best result we knew was of the form, see \cite{ch1}: if $m=\frac{n}{2}+2$ and $s=\frac{n}{2}+4$ then the quadratic (D-NLS) is locally well-posed.  This is a bit too far from the scaling exponent and, as we will see later on, the decay is too strong also.

The analysis of the problem brings the conclusion that the ``worst'' interactions are the orthogonal ones, i.e. those between waves which travel in orthogonal directions. Therefore the problem becomes more interesting in dimensions $2$ or higher and we decided to understand what happens when $n=2$. This is why in this work we decided to specialize to the case of two-dimension quadratic (D-NLS).   
 
Our goal in the first place was to obtain local well-posedness for initial data $u_{0} \in H^{s}$, for any $s > s_{c}$. To achieve that goal we assumed that the initial data comes with a bit of spherical symmetry.

We recall the definition of the differential operator: 
\beq \label{r}
\mathcal{R}f=(x_{1} \partial_{x_{2}} - x_{2} \partial_{x_{1}})f 
\eeq

\noindent
and of the pseudo differential operators in the left calculus:

\beq \label{de}
\mathcal{D}'f= (1+\frac{\langle x \rangle^{2}}{\mu+\langle D \rangle^{2}})^{\frac{1}{4}+\frac{\e}{2}}f \ \ \mbox{with symbol} \ \ (1+\frac{\langle x \rangle^{2}}{\mu+\langle(\xi,\tau) \rangle^{2}})^{\frac{1}{4}+\frac{\e}{2}}
\eeq

\noindent 
for some $0 < \e < \q$. For a generic space of functions $\mathcal{X}$ we defined:

\beq \label{S}
\mathcal{D}'\mathcal{R}\mathcal{X}= \{f \in \mathcal{X}: \mathcal{D}'f\in \mathcal{X} \ \mbox{and} \  \mathcal{D}'\mathcal{R}f \in \mathcal{X}\}
\eeq

We renamed the decay operator from Part 1 by $\mathcal{D}'$ since we use in the present paper a more general type of decay which we call $\mathcal{D}$. 

The main result of Part 1 is the following:

\begin{t1} \label{TT}

Assume n=2. Given any $s > s_{c}$ and $T > 0$, there exists $\delta > 0$ such that for every $u_{0} \in \mathcal{D}'\mathcal{R}H^{s}$ with $\delta_{0}=||u_{0}||_{\mathcal{D}'\mathcal{R}H^{s}}  < \delta$ , the quadratic (D-NLS) has a unique solution $u$ in $C([0,T]:\mathcal{D}'\mathcal{R}H^{s}) \cap \mathcal{R}\mathcal{D}' Z^{s,5}$ with Lipschitz dependence on the initial data.

\end{t1}

The definition of $Z^{s,5}$ will come up in the current paper. Two major questions arise once we acknowledge this result. One is to try to obtain a result without involving any spherical symmetry and the other one is to remove the smallness condition on the initial data.

The current paper answers to the first issue. We essentially prove that for any $s > s_{c}+1$ the quadratic (D-NLS) is locally well-posed for small $u_{0} \in \mathcal{D}H^{s}$. In section \ref{ds} we provide a precise definition of $\mathcal{D}H^{s}$; if the reader digested the definition of $\mathcal{D}'$ we can remark that $\mathcal{D}'H^{s} \subset \mathcal{D}H^{s}$. 

Let us make the result we obtain more precise.

We denote  by $\chi_{[0,T]}$ a smooth approximation of the characteristic function of $[0,T]$ such that $\chi_{[0,T]}(t)=1, \ \forall t \in [0,T]$. We will always consider $\chi_{[0,T]}$ as a function of time, in other words by $\chi_{[0,T]}$ we mean  $\chi_{[0,T]}(t)$.

We dedicate the section \ref{ds} to the definition of the spaces  $\mathcal{D}Z^{s,5}$ (for the solutions) and $\mathcal{D}W^{s}$ (for the inhomogeneity). These spaces satisfy the linear estimate:

\begin{t1} \label{le1}

If $g \in \mathcal{D}H^{s}$ and $f \in \mathcal{D}W^{s}$, then the solution of:

\begin{equation} \label{E}
\begin{cases}
\begin{aligned} 
&iu_{t}-\Delta u = f \\
&u(x,0) =g(x)
\end{aligned}
\end{cases}
\end{equation}

satisfies  $\chi_{[0,1]} u \in \mathcal{D}Z^{s,5} \cap C_{t}\mathcal{D}H_{x}^{s}$.

\end{t1}

To each quadratic nonlinearity we associate is the standard way the bilinear form $B_{P}(u,v)$. The bilinear estimate is the next key result:

\begin{t1} \label{bg}
If $s > s_{c}+1$, we have the global bilinear estimate:

\beq \label{bg1}
||B_{P}(u,v)||_{\mathcal{D}W^{s}} \leq C_{s} ||u||_{\mathcal{D}Z^{s}}||v||_{\mathcal{D}Z^{s}} 
\eeq

\end{t1}

Once we have the above two results, a standard fixed point argument gives us the main result:

\begin{t1} \label{T}

Assume n=2. Given any $s > s_{c}+1$ and $T > 0$, there exists $\delta > 0$ such that for every $u_{0} \in \mathcal{D}H^{s}$ with $\delta_{0}=||u_{0}||_{\mathcal{D}H^{s}}  < \delta$ , the quadratic (D-NLS) has a unique solution $u$ in $C([0,T]:\mathcal{D}H^{s}) \cap \mathcal{D} Z^{s,5}$ with Lipschitz dependence on the initial data.

\end{t1}

The general approach of this result is similar to the one in Part 1. Let $B(u,v)$ be the bilinear form:

\beq \label{bf}
B(u,v)=\sum_{i,j \in \{1,2\}} c_{ij} u_{x_{i}} v_{x_{j}}
\eeq

\noindent
where $c_{ij}$ are complex constants. We intend to obtain bilinear estimates for $B(u,v)$ and $B(u,\bar{v})$ since this way we cover the theory for all quadratic polynomials of type $P(\nabla{u},\nabla{\bar{u}})$, except for those of type $P(\nabla{\bar{u}})$. For the last ones the theory had been developed previously, see \cite{g1}.

 We start with $X^{s,\q,1}$ as the candidate for $Z^{s}$ and $X^{s,-\q,1}$ as a candidate for $W^{s}$. The bilinear estimates work fine as long as we recover information which is at some distance from the paraboloid ($\tau=\xi^{2}$) and it breaks down very close to the paraboloid - we catch a logarithm of the high frequency which cannot be controlled. To remedy this we come up with a more delicate decomposition of the part of the Fourier space which is at distance less than $1$ from the paraboloid. More exactly we introduce a wave packet decomposition and we measure the packets in $L^{\infty}_{t}L^{2}_{x}$. Then the target space $W^{s}$ is also modified at distance less than $1$ from paraboloid, i.e. we also have a wave packet decomposition and the packets are measured in $L^{1}_{t}L^{2}_{x}$. We have to recover a $L^{1}_{t}$ structure on the packets for $B(u,v)$ and this is why we need to involve the extra decay.

All along the argument we do involve decay in the bilinear estimates and this is why our spaces will be of type $\mathcal{D}Z^{s}$ and $\mathcal{D}W^{s}$. See section \ref{ds} for the definitions.

Once the bilinear estimates are fixed, then a standard fixed point argument gives us the result of Theorem \ref{T}. 

One can easily adapt our argument for the bilinear forms of type:

\beq \label{bf1}
B(u,v)=\sum_{j=1}^{2} c_{j} u v_{x_{j}}
\eeq

 This is because the basic estimates are derived for the bilinear form $\tilde{B}(u,v)=u \cdot v$ and then we "over-estimate" the size of $\nabla$, see the beginning of section \ref{bes} for more details. Thus we are entitled to claim the result for the quadratic polynomials of type $P(u,\nabla{u})$,  $P(\bar{u},\nabla{u})$ and $P(u,\nabla{\bar{u}})$.

The spaces we use in this paper are in some way the counterpart of the ones involved in dealing with the wave maps equation, see \cite{ta1} and \cite{tao1}. Our spaces are a bit more difficult since they involve phase-space localization, rather than phase localization which is the case for the wave-maps. 

We should also make a point in the fact that the result in the current paper is not a trivial reproduction of the argument in Part 1, for the case when we do not use any spherical symmetry. One would notice along the proof that we need to use decay when obtaining estimates solely in $X^{s,\q,1}$ spaces. In Part 1, we were able to derive the bilinear estimates in $X^{s,\q,1}$ by using only the spherical symmetry. We also changed the type of decay and not just to make it more general. In Part 1 we really needed an hypoelliptic operator to give decay, while in the current paper the decay we use some sort of a micro-local version of the condition in  (\ref{co1}). The current decay type is more general than the one used in Part 1; on the other hand it would not be good enough for the purpose there.

We conclude the introduction with few open problems. We predicted from Part 1 that, without assuming any symmetry, we do expect a positive result for $s > s_{c}+1$ and a negative one for $s < s_{c}+1$. In this paper we provide only the positive result; the negative one is subject to current research. 

The generalization to higher dimensions should be of interest too. We know that the scaling exponent is $\frac{n}{2}$ for the case when both terms come with derivatives and we think it should be possible to get similar results under similar conditions in all dimensions.

\vspace{.2in}

\section{Definition of the spaces} \label{ds}

For each $u$ we denote by $\mathcal{F}u=\hat{u}$ the Fourier transform of $u$. This is always taken with respect to all the variables, unless otherwise specified.

Throughout the paper $A \lesssim B$ means $A \leq CB$ for some constant $C$ which is independent of any possible variable in our problem. We say $A \approx B$ if $A \leq CB \leq C^{2}A$ for the same constant $C$. We say that we localize at frequency $2^{i}$ to mean that  in the support of the localized function $|(\xi,\tau)| \in [2^{i-1},2^{i+1}]$.

In the Schr\"odinger equation time and space scale in a different way, and this suggests to define the norm for $(\xi, \tau)$ by $|(\xi, \tau)|=(|\tau|+ \xi^{2})^{\q}$. In dealing with the quadratic nonlinearity without derivatives the Bourgain space $X^{s,b}$ proved to be a very useful. They are defined in the following way:

$$X^{s,b} =  \{f\in S'; \langle (\xi,\tau) \rangle^{s} \langle \tau-\xi^{2} \rangle^{b} \hat{f} \in L^{2} \}$$

Here and thereafter $\langle x \rangle=(1+|x|^{2})^{\q}$ where $|x|$ is the norm of $x$. We will employ frequency localized versions of $X^{s,\q}$ which are constructed according to weights present in its definition. 

Consider $\varphi_{0} : [0,\infty) \rightarrow \R$ be a nonnegative smooth function such that $\varphi_{0}(x)=1$ on $[0,1]$ and $\varphi_{0}(x)=0$ if $x \geq 2$. Then for each $i \geq 1$ we define $\varphi_{i}: [0,\infty) \rightarrow \R$ by $\varphi_{i}(x)= \varphi_{0}(2^{-i}x)- \varphi_{0}(2^{-i+1}x)$. We define the operators $S_{i}$, to localize at frequency $2^{i}$, by:

$$\mathcal{F}(S_{i}f)= \hat{f}_{i} = \varphi_{i}(|(\xi, \tau)|) \cdot \hat{f}(\xi, \tau)$$

For $d \in I_{i}=\{2^{-i}, 2^{-i+1}, .., 2^{i+2}\}$ we define $\varphi_{i,d}(\xi,\tau)=\varphi_{i}(|(\xi, \tau)|) \cdot \varphi_{i+\ln_{2}{d}}(|\tau-\xi^{2}|)$. There is one simple reason to chose to work with $d$ in this way rather than working with $2^{d}$. If $|(\xi,\tau)| \approx 2^{i}$ then $|\tau-\xi^{2}| \approx |(\tau,\xi)| d((\xi,\tau),P) \approx 2^{i} d((\xi,\tau),P)$ (away from zero). Hence one should think of $d$ as the distance to $P$ since the support of $\varphi_{i,d}$ is approximately the set 

$$\{(\xi,\tau): |(\xi,\tau)| \approx 2^{i}, d((\tau,\xi),P) \approx d\} \approx \{(\xi,\tau): |(\xi,\tau)| \approx 2^{i}, |\tau -\xi^{2}| \approx d2^{i}\}$$

It is easy to notice that

$$ \sum_{d \in I_{i}} \varphi_{i,d} (\xi,\tau) = \varphi_{i}((|\xi,\tau)|)  , \ \forall (\xi,\tau) \in \R^{2} \times \R $$

We define the operators $S_{i,d}$ by $S_{i,d}f=f_{i,d}=\check{\varphi}_{i,d} * S_{i}f$ and we have $f_{i} = \sum_{d \in I_{i}} f_{i,d}$. In the support of $\hat{f}_{i,d}$ we have $1+|\tau - \xi^{2}| \approx 2^{i}d$.

Sometimes it is useful to localize in a linear way rather than a dyadic way. In these cases we localize with respect to the value of $|\tau - \xi^{2}|$ instead; we will make this clear when we need it.

For each dyadic value $d \in I_{i}$ we introduce the operators which localize at distance less and greater than $d$ from $P$: 

$$S_{i, \leq d}f = f_{i, \leq d}=\sum_{d' \in I_{i}: d' \leq d} f_{i, d'} \ \  \mbox{and} \ \ S_{i, \geq d}f=f_{i, \geq d}= f_{i}-f_{i, \leq d}$$

The part of $\hat{f}$ which is at distance less than $1$ from $P$ plays an important role and this is why we define the global operators:

$$S_{\cdot,\leq 1}f=f_{\cdot, \leq 1}=\sum_{i=0}^{\infty} f_{i, \leq 1} \ \  \mbox{and} \ \ S_{\cdot,\geq 1}f=f_{\cdot,\geq 1}=\sum_{i=0}^{\infty} f_{i, \geq 1}$$

We denote by $A_{i}$ the support in $\R^{2} \times \R$ of $\varphi_{i}(|(\xi,\tau)|)$ and by $A_{i,d}$ the support of $\varphi_{i,d}$. In a similar way we can define $A_{i, \leq d}$ and $A_{i, \geq d}$ to be the support of the operators $S_{i, \leq d}$, respectively $S_{i, \geq d}$.

For functions whose Fourier transform is supported in $A_{i}$ we define and for any $1 \leq p \leq \infty$: 

$$||f||^{p}_{X_{i}^{0,\q,p}} = \sum_{d \in I_{i}} ||S_{i,d}f||_{X^{0,\q}}^{p}$$

\noindent
with the usual convention for $p=\infty$. Then we define the space $X^{s,\q,p}$ by the norm:

$$||f||^{2}_{X^{s,\q,p}}=\sum_{i} 2^{2is} ||f_{i}||_{X^{0,\q,p}_{i}}$$

For technical purposes we need localized versions of this spaces, like $ X^{s,\q}_{i,d} = \{f \in X^{s,\q}: \hat{f} \ \mbox{supported in} \ A_{i,d} \}$ and, similarly, $ X^{s,\q,p}_{i,\leq d}$ and $ X^{s,\q,p}_{i,\geq d}$.

$X^{s,\q,1}$ is our first candidate for the space of solutions. Our computations indicate that it is the right space to measure only the part of the solution whose support in the Fourier space is at distance greater than $1$ from $P$, i.e. the $S_{\cdot,\geq 1}$ part of our solutions. 
 
The $S_{\cdot,\leq 1}$ part of the solutions can be measured in $X^{s,\q,\infty}$ plus an additional structure whose construction is described bellow. 

We define the following lattice in the plane $\tau=0$:

$$\Xi = \{\xi=(r,\t): r=n,\  \t=\frac{\pi}{2} \frac{k}{n}, \ n,k \ \mbox{positive integers} \}$$

$\Xi$ is like a lattice in polar coordinates. It has the properties that the distance between any two points is at least $1$ and that for every $\eta \in \R^{2}$ there is a $\xi \in \Xi$ such that $|\xi-\eta| \leq 1$. For each $\xi \in \Xi$ we build a non-negative function $\phi_{\xi}$ to be a smooth approximation of the characteristic function of the cube of size $1$ in $\R^{2}$ centered at $\xi$ and satisfying the natural partition property:

$$\sum_{\xi \in \Xi} \phi_{\xi}=1$$

We can easily impose uniforms bounds on the derivatives of the system $(\phi_{\xi})_{\xi \in \Xi}$. For each $\xi \in \Xi$ we define:

$$f_{\xi} = \check{\phi}_{\xi} * f \ \ \ \ \mbox{and} \ \ \ f_{\xi, \leq 1} = \check{\phi}_{\xi} * f_{\cdot,\leq 1}$$

The convolution above is performed with respect to the $x$ variable, i.e. it does not involve the $t$ variable. The support of $\hat{f}_{\xi, \leq 1}$ is like a parallelepiped having the center $(\xi,\xi^{2}) \in P$ and sizes: $\approx |\xi|$ in the $\tau$ direction and $1$ in the other two directions (normal to $P$ and the completing third one).

The next concern is how to measure $f_{\xi, \leq 1}$. Let's denote by $(Q^{m})_{m \in \Z^{2}}$ the standard partition of $\R^{2}$ in cubes of size $1$; i.e. $Q^{m}$ is centered at $m=(m_{1},m_{2}) \in \Z^{2}$, has its sides parallel to the standard coordinate axis and has size $1$. For each $\xi \in \Xi$, $m \in \Z^{2}$ and $l \in \Z$ we define the tubes:

$$T_{\xi}^{m,l}=\cup_{t \in [l,l+1]} (Q^{m}-2t\xi) \times \{ t \}=$$

$$\{(x-2t\xi_{1}, y-2t\xi_{2},t): (x,y) \in Q^{m} \ \mbox{and} \ t \in [l,l+1] \}$$

Then, for each $\xi \in \Xi$, we define the space $Y_{\xi}$ by the following norm:

$$||f||^{2}_{Y_{\xi}}= \sum_{(m,l) \in \Z^{3}} ||f||^{2}_{L^{\infty}_{t}L^{2}_{x}(T_{\xi}^{m,l})}$$

We have $f=\sum_{\xi \in \Xi} f_{\xi}$ and then we define the space $Y^{s}$ by the norm:

$$||f||^{2}_{Y^{s}}= \sum_{\xi \in \Xi} \langle \xi \rangle^{2s} ||f_{\xi}||^{2}_{Y_{\xi}}$$

We define also the localized versions $Y_{i}= \{ f \in Y^{0}; \hat{f} \ \mbox{supported in} \ A_{i} \}$ and $Y_{i, \leq d}= \{ f \in Y^{0}; \hat{f} \ \mbox{supported in} \ A_{i, \leq d} \}$, the last one being defined for any $d \in I_{i}$ with $d \leq 1$.

Our solutions will be localized in time. If we come with a frequency localization on the top of this we are left with decay in time of our solutions. For this we define $Y_{\xi}^{N}$ and $Y^{s,N}$ by the norms:

$$||f||_{Y_{\xi}^{N}}=||\langle t \rangle^{N}f||_{Y_{\xi}} \ \ \ \mbox{and} \ \ \ ||f||^{2}_{Y^{s,N}}= \sum_{\xi \in \Xi} \langle \xi \rangle^{2s} ||f_{\xi}||^{2}_{Y_{\xi}^{N}}$$

To bring everything together, define $Z^{s,N}$ to be 

$$Z^{s,N}=\{ f \in S': ||f_{\cdot,\geq 1}||_{X^{s,\q,1}}+||f_{\cdot, \leq 1}||_{Y^{s,N}}+||f_{\cdot, \leq 1}||_{X^{s,\q,\infty}} < \infty \}$$

\noindent
with the obvious norm. One important property, proved in Part 1, is:

$$X^{s,\q,1} \subset Z^{s}$$

Our spaces are equipped with an additional decay structure which we describe bellow. For each $i$, let $Q^{m}_{i}$ be a system of cubes of size $2^{i}$ which form a partition of $\R^{2}$; we choose them so that their sides are parallel to the coordinate axes and the center of $Q^{m}_{i}$ is $(2^{i}m_{1},2^{i}m_{2})$. Let $L=\{(x,y):ax+by=0\}$ be the equation of a line passing through the origin and denote by $\bar{n}$ the normal unit vector to $L$. For each $k \in \Z$, we define $L^{k}_{i}=\{(x',y'): (x',y')=(x,y)+k 2^{i} \bar{n}; (x,y) \in L\}$ to be the line parallel to $L$ and at distance $2^{i}k$ from $L$. If $f : \R^{2} \to \C$ we introduce the norm:

$$||f||_{\mathcal{D}_{i}L^{2}}=\sup_{L} \sup_{k \in \Z} \left( \sum_{m: Q_{i}^{m} \cap L^{k}_{i} \ne \emptyset} ||f||_{L^{2}(Q_{i}^{m})} \right)$$

What happen above is that we sum up the $||f||_{L^{2}(Q_{i}^{m})}$ in $l^{1}$ over those $Q_{i}^{m}$'s which intersect $L^{k}_{i}$, and then we take a suppremum with respect to $k$. In the end we take a suppremum with respect to all lines $L$.

Our decay space is defined by the norm:

$$||f||_{\mathcal{D}H^{s}}^{2}=\sum_{i} 2^{2is} (||f_{i}||^{2}_{\mathcal{D}_{i}L^{2}} + ||f_{i}||^{2}_{L^{2}})$$

If $f: \R^{3} \to \C$ (we include the time dependent functions), we define:

$$||f||_{\mathcal{D}_{i}L^{2}}=\sup_{L} \sup_{k \in \Z} \left( \sum_{m: Q_{i}^{m} \cap L^{k}_{i} \ne \emptyset} ||f||_{L^{2}(Q_{i}^{m} \times \R)} \right)$$

For $f$ such that $\hat{f}$ is supported in $A_{i,d}$ we define:

$$||f||_{\mathcal{D}_{i}X^{s,\q}_{i,d}}=2^{js} (2^{j}d)^{\q} (||f||_{\mathcal{D}_{i}L^{2}}+||f||_{L^{2}})$$

\noindent
and the decay version of $X^{s,\q,p}$,  $\mathcal{D}X^{s,\q, p}$ is defined by the norm:

$$||f||_{\mathcal{D}X^{s,\q, 1}}^{2}=  \sum_{i} (\sum_{d} ||f_{i,d}||^{p}_{\mathcal{D}_{i}X^{s,\q}_{i,d}})^{\frac{2}{p}}$$

 For $f \in Y_{i}$ we define:

$$||f||_{\mathcal{D}_{i}Y_{i}}=\sup_{L} \sup_{k \in \Z} \left( \sum_{m: Q_{i}^{m} \cap L^{k}_{i} \ne \emptyset} (\sum_{\xi: |\xi| \approx 2^{i}} ||\chi_{Q_{i}^{m} \times \R}f_{\xi}||^{2}_{Y_{\xi}})^{\q} \right)+||f||_{Y}$$

The decay version of $Y^{s}$, $\mathcal{D}Y^{s}$ is defined by the norm:

$$||f||^{2}_{\mathcal{D}Y^{s}}= \sum_{i} 2^{2is} ||f_{i}||^{2}_{\mathcal{D}_{i}Y_{i}}$$

To bring everything together, define $\mathcal{D}Z^{s,N}$ to be 

$$\mathcal{D}Z^{s,N}=\{ f \in S': ||f_{\cdot,\geq 1}||_{\mathcal{D}X^{s,\q,1}}+||f_{\cdot, \leq 1}||_{\mathcal{D}Y^{s,N}}+||f_{\cdot, \leq 1}||_{\mathcal{D}X^{s,\q,\infty}} < \infty \}$$

So far we have built the spaces suitable for the solution of (\ref{E:s}). We need also a space for the right hand side of the equation, see Theorem \ref{le1}.

We can easily define $X^{s,-\q,p}$ by simply replacing $\q$ with $-\q$ in the definition of $X^{s,\q,p}$. Then we define $\mathcal{Y}^{s}$ by:

$$||f||^{2}_{\mathcal{Y}^{s,N}}= \sum_{\xi \in \Xi} \langle \xi \rangle ^{2s} ||f_{\xi}||^{2}_{\mathcal{Y}_{\xi}^{N}}$$

\noindent
where $\mathcal{Y}_{\xi}^{N}$ is defined as follows:

$$||f||^{2}_{\mathcal{Y}_{\xi}^{N}}= \sum_{(m,l) \in Z^{3}} || \langle t \rangle^{N} f||^{2}_{L^{1}_{t}L^{2}_{x}(T_{\xi}^{m,l})}$$

Notice that $(\mathcal{Y}_{\xi})^{*}=Y_{\xi}$ since we will use this later for duality purposes.

 We introduce $W^{s,N}$ defined by the norm:

$$||f||_{\mathcal{W}^{s,N}}= \inf \{||f_{1}||_{\mathcal{Y}^{s,N}} + ||f_{2}||_{X^{s,-\q,1}}; f=f_{1}+f_{2} \}$$

We measure the right hand side of ($\ref{E}$) in: 

$$W^{s,5} =\{ f \in S': ||f_{\cdot,\leq 1}||^{2}_{\mathcal{W}^{s,5}} + ||f_{\cdot,\geq 1}||^{2}_{X^{s,-\q,1}} < \infty \}$$

As before, we can define $\mathcal{D}X^{s,-\q,1}$, $\mathcal{D}\mathcal{Y}^{s,N}$ and $\mathcal{D}W^{s,5}$. 

Besides $X^{s,b}$ we need the conjugate $\bar{X}^{s,b}$ which is defined as follows:

$$\bar{X}^{s,b} =  \{f\in S'; \langle (\xi,\tau) \rangle^{s} \langle \tau+\xi^{2} \rangle^{b} \hat{f} \in L^{2}  \}$$

We can define all the other elements the same way as above by simply placing a bar on each space and operator, while replacing everywhere $|\tau-\xi^{2}|$ with $|\tau+\xi^{2}|$ and $P$ with $\bar{P}=\{(\xi,\tau): \tau+\xi^{2}=0 \}$.  

We record the following important facts:

$$ f \in X^{s,b} \Longleftrightarrow \bar{f} \in \bar{X}^{s,b} \ \ \ \mbox{and} \ \ \ (X^{0,\q})^{*}=\bar{X}^{0,-\q}$$

Before we start we need to introduce some new localization operators. For each $i \in \N$ we define a refined lattice:

\beq \label{Xi}
\Xi^{i}=\{\xi=(n2^{-i},\t): \t=\frac{\pi}{2} \frac{l}{n}; n \in \N, l \in \Z \}
\eeq

For each $\xi \in \Xi^{i}$ we build  the corresponding $\phi^{i}_{\xi}$ to be a smooth approximation of the characteristic function of the cube centered at $\xi$ and with sizes $2^{-i}$. We also assume that the system $(\phi_{\xi}^{i})_{\xi \in \Xi^{i}}$ forms a partition of unity in $\R^{2}$.

For each $l \in \Z$ we can easily construct a function $\chi_{[l-\q,l+\q]}$ to be a smooth approximation of the characteristic function of the interval $[l-\q,l+\q]$ and such that the system $(\chi_{[l-\q,l+\q]})_{l \in \Z}$ form a partition of unity in $\R$. For any $\xi \in \Xi^{i}$ with $|\xi| \leq 2^{i+1}$ we consider those $l \in \Z$ with the property $|(\xi,l)| \approx 2^{i}$ and define the operators: 

$$\hat{f}_{\xi,l}= \phi_{[l-\q,l+\q]}(\tau) \phi^{i}_{\xi}(\xi) \hat{f}(\xi,\tau) $$  

  The support of $\hat{f}_{\xi, l}$ is approximately a tube centered at $(\xi,l)$ and of size $2^{-i} \times 2^{-i} \times 1$, the last one being in the $\tau$ direction. Since the distance of these tubes will play an important role, sometimes it would be convenient if we were able to work with $(f_{\xi, \xi^{2}+l})_{\xi \in \Xi^{i}, l \in \Z}$ instead. The only problem is that it is not guaranteed that $\xi^{2} \in \Z$ for all $\xi \in \Xi^{i}$. Of course we could change the way we cut in the $\tau$ direction, but this would complicate notations even more. We choose instead to ignore that $\xi^{2}$ may not be integer, and go on and use $g_{\xi, \xi^{2}+l}$. It will be obvious from the argument that this does not affect in any way the rigorousness of the proof. The last notation we introduce is $f_{\xi, \xi^{2} \pm l}=f_{\xi, \xi^{2} + l}+f_{\xi, \xi^{2} - l}$.

For $d \leq 2^{i-2}$ we obtain a new decomposition of $g_{i,d}$:

\beq \label{aux2}
g_{i, d} = \sum_{k =2^{i-1}d}^{2^{i+1}d} \sum_{\xi \in \Xi^{i}} g_{\xi, \xi^{2} \pm k}
\eeq

Notice that the $\xi$'s $\in \Xi^{i}$ involved in the above summation have $|\xi| \approx 2^{i}$.

For the part of $\hat{g}$ supported away from $P$ we come with a different decomposition: 

\beq \label{aux1}
g_{i, \geq 2^{i-2}} = \sum_{n} \sum_{\xi \in \Xi^{i}} \sum_{l \in I_{\xi}} g_{\xi,l}
\eeq

\noindent
where $I_{\xi}=\{l \in \Z: 2^{2i-2} \leq |l-\xi^{2}| \leq 2^{2i+2} \}$. The $\xi$'s $\in \Xi^{i}$ involved in the above summation have $|\xi| \leq 2^{i+1}$.

\section{Proof of Theorem \ref{T} and Linear Estimates}

In this section we intend to use the results already proved in Part 1, see \cite{be}. It is well-known in the literature that once we have the linear estimates in Theorem \ref{le1} and the bilinear estimates in Theorem \ref{bg} we get the result in Theorem \ref{T} via a standard fixed point argument. See, for instance, Part 1.

The linear estimates, see Theorem \ref{le1}, in the variant without decay, were also proved in Part 1, see the corresponding section there. The conservation of decay can be easily adapted to those proofs. This is due to the fact that the type of decay we use is scaled properly for the Schr\"odinger equation. If one wants to pursue the complete argument, it would be useful to acknowledge the section results in section \ref{ar}.

\section{Bilinear estimates in $\mathcal{R}X^{s,\q,1}$} \label{bes}

In this section we derive the bilinear estimates for $B(u,v)$ and $B(u,\bar{v})$ in $\mathcal{R}X^{s,\q,1}$, where $B$ is of type (\ref{bf}). We introduce the additional bilinear form $\tilde{B}(u,v)=u \cdot v$. If $\hat{u}$ is localized in $A_{i}$ we use the estimate $||\nabla{u}||_{L^{2}} \leq 2^{i} ||u||_{L^{2}}$. $X^{s,\pm \q,1}$ are $L^{2}$ like on dyadic pieces, hence if $\hat{u}$ is localized in $A_{i}$ and $\hat{v}$ is localized in $A_{j}$ we use the estimates:

\beq \label{bb}
 ||\tilde{B}(u,v)||_{X} \leq C ||u||_{X'} ||v||_{X''}   \Ra  ||B(u,v)||_{X} \leq 2^{i+j} C ||u||_{X'} ||v||_{X''}
\eeq

\beq \label{bbb}
 ||\tilde{B}(u,\bar{v})||_{X} \leq C ||u||_{X'} ||v||_{X''}   \Ra  ||B(u,\bar{v})||_{X} \leq 2^{i+j} C ||u||_{X'} ||v||_{X''}
\eeq

Here $X,X',X''$ are of type $X^{s,\pm \q, p}$. The constant $C$ may depend on $u,v$, more exactly of their localizations. The key thing is once we have estimates for $\tilde{B}$, we obtain estimates for $B$ by simply bringing in the correction factor of $2^{i+j}$.  

If $B$ were of type (\ref{bf1}) the correction factor would be only $2^{j}$ and this justifies why we can claim the estimates for bilinear estimates of type (\ref{bf1}).

Another thing to keep in mind is that we apply duality along the proof and anytime we do it we mean it in the context of  $\tilde{B}$, not $B$.

The main results we claim are listed in the following theorem.

\begin{t1} \label{bil}

a) If $i \leq j$, we have the following estimates:

\beq \label{be11}
||B(u,v)||_{\mathcal{D}_{k}X^{s,-\q}_{k}} \les j 2^{(2-s)i} 2^{(k-j)s}||u||_{\mathcal{D}_{i}X^{s,\q}_{i}} ||v||_{\mathcal{D}_{j}X^{s,\q}_{j}}
\eeq

The above estimates holds true if $B(u,v)$ is replaced by $B(\bar{u},v)$ or $B(u,\bar{v})$. 

b) If $5i \leq j$, we have the following estimates:

\beq \label{be22}
||B(u,v)||_{\mathcal{D}_{j}X^{s,-\q}_{j, \geq 2^{-i}}} \les 2^{(2-s)i} i||u||_{\mathcal{D}_{i}X^{s,\q}_{i}} ||v||_{\mathcal{D}_{j}X^{s,\q}_{j, \geq 2^{-i}}}
\eeq

The above estimates holds true if $B(u,v)$ is replaced by $B(\bar{u},v)$ or $B(u,\bar{v})$.

\end{t1}

\vspace{.2in}

\subsection{Basic Estimates}
\noindent

\vspace{.1in}

We start with a simple result stating how two parabolas interact under convolution. We need few technical definitions. 

Throughout this section functions are defined on Fourier space (they should be thought as Fourier transforms). This is why we use the standard coordinates $(\xi,\tau)$. 

For each $c \in \R$ denote by $P_{c}= \{(\xi,\tau): \tau-\xi^{2}=c\}$ and  by $\bar{P}_{c}=\{(\xi,\tau): \tau+\xi^{2}=c \}$. For simplicity $P=P_{0}$ and $\bar{P}=\bar{P}_{0}$. 

Denote by $\delta_{P_{c}}=\delta_{\tau-\xi^{2}=c}$ the standard surface measure associated to the parabola $P_{c}$. With respect to this measure, the restriction of $f$ to $P_{c}$ has norm:

$$||f||_{L^{2}(P_{c})}=\left( \int f^{2}(\xi,\xi^{2}+c) \sqrt{1+4|\xi|^{2}} d\xi  \right)^{\q}$$

The first result was derived in Part 1, see the corresponding section there:

\begin{p3} \label{pp1}

 Let $f \in L^{2}(P^1)$ and  $g \in L^{2}(P^2)$ such that $f$ is localized at frequency $2^{i}$ and $g$ at frequency $2^{j}$. We have:

\beq \label{ge}
||f \delta_{P^{1}} * g \delta_{P^{2}} ||_{L^{2}} \les 2^{\min{(i,j)}} ||f||_{L^{2}(P^{1})}||g||_{L^{2}(P^{2})}
\eeq
 
\noindent
where $P^{1} \in \{P_{c_{1}},\bar{P}_{c_{1}}\}$ and $P^{2} \in \{P_{c_{2}},\bar{P}_{c_{2}}\}$.

\end{p3}

The second result comes to replace the corresponding one in Part 1 for the case when we do not have any symmetry involved.

\begin{p3} \label{pp2}

 We assume that we are in the same setup as in Proposition \ref{pp1}. In addition we assume $i \leq j$, $|c_{1}| \leq 2^{2i-2}$ and $|c_{2}| \leq 2^{i+j-4}$. Then 

\beq \label{ge1}
||f \delta_{P^{1}} * g \delta_{P^{2}} ||_{L^{2}(|(\xi,\tau)| \approx 2^{j},|\tau-\xi^{2}| \leq d)} \les d^{\q} ||f||_{L^{2}(P^{1})}||g||_{L^{2}(P^{2})}
\eeq
 
\noindent
where $P^{1} \in \{P_{c_{1}},\bar{P}_{c_{1}}\}$ and $P^{2} \in \{P_{c_{2}},\bar{P}_{c_{2}}\}$.

\end{p3}

\begin{proof}

We notice that it is enough to prove the result under the hypothesis that $d \leq 2^{i+j-4}$ since otherwise, the result in (\ref{ge}) is stronger. Without losing generality we can assume $c_{1}=c_{2}=0$. One could easily adapt the argument bellow to the general case when $|c_1|+|c_2| \leq 2^{i+j-2}$. 

{\mathversion{bold} $f \delta_{P} * g \delta_{P}$}

An easy way to test the above norm of a function is by estimating $|(f \delta_{\tau=\xi^{2}} * g \delta_{\tau=\xi^{2}})h|$ for any $h \in L^{2}$ supported in the region $|\tau-\xi^{2}| \leq d$. For any such $h$ we have:

$$(f \delta_{\tau=\xi^{2}} * g \delta_{\tau=\xi^{2}})h=$$

$$\int f(\xi) g(\eta) h(\xi+\eta, \xi^{2}+\eta^{2}) \sqrt{1+4\xi^{2}}  \sqrt{1+4\eta^{2}} d\xi d\eta $$

Since $h$ is supported in a region $|\tau-\xi^{2}| \leq d$ we need the following condition on the variables inside the integral: $|(\xi+\eta)^{2}-(\xi^{2}+\eta^{2})| \leq d$ or $2|\xi||\eta|\cos{\t} \leq d$ where $\t$ is the angle between $\xi$ and $\eta$. 
Hence $\cos{\t} \leq 2^{-i-j}d$ which implies $|\t-\frac{\pi}{2}| \leq 2^{-i-j}d$. This suggests decomposing:

$$[0,2\pi]=\bigcup_{l \in J_{i,j,d}} I_{l}=\bigcup_{l=1}^{2^{i+j+2}d^{-1}} [(l-\q)2^{-i-j}d\frac{\pi}{2},(l+\q)2^{-i-j}d \frac{\pi}{2}]$$

\noindent
in other words to split $[0,2\pi]$ in a disjoint union of intervals of size $2^{-i-j}d \frac{\pi}{2}$. Correspondingly we split:

$$f=\sum_{l \in J_{i,j,d}} f_{l}  \ \ \ \mbox{and} \ \ \ g=\sum_{l \in J_{i,j,d}} g_{l}$$

\noindent
such that $f_{l}$ is the part of $f$ localized in $A_{l}=\{\xi: \arg{\xi} \in I_{l} \}$ and similarly for $g$.

If $\arg \xi \in I_{l}$ and $\arg \eta \in I_{l'}$ and we want them to belong to the domain of integration above we need to impose $|l-l'|=2^{i+j}d^{-1}$ (modulo $2^{i+j+2}d^{-1}$). For each $l$ there are two $l'$'s with this property. We simplify more and consider that there is only one $l^{\perp}=l+2^{i+j}d^{-1}$ (modulo $2^{i+j+2}d^{-1}$) with this property; one can easily complete the argument with both values. Then we have:

$$(f \delta_{\tau=\xi^{2}} * g \delta_{\tau=\xi^{2}})h=$$

$$\sum_{l \in J_{i,j,d}} \int f_{l}(\xi) g_{l^{\perp}}(\eta) h(\xi+\eta, \xi^{2}+\eta^{2}) \sqrt{1+4\xi^{2}}  \sqrt{1+4\eta^{2}} d\xi d\eta $$

Now that we have a sharp angular localization, we complete it with a norm localization which should be consistent with the angular one:

$$f_{l}=\sum_{m} f_{l,m} \ \ \ \mbox{and} \ \ \ g_{l^{\perp}}=\sum_{n} g_{l^{\perp},n}$$

One can easily notice that it is not important to explicitly delimit the sets for $m$ and $n$.

For the low frequency things are simple: $f_{l,m}$ is the part of $f_{l}$ localized in the set  $A^{i}_{l,m}=\{\xi \in A_{l}:|\xi| \in [2^{-i}d(m-\q),2^{-i}d(m+\q)] \}$ and notice that this is consistent with the arc length size localization of $g_{l^{\perp}}$ (which is $2^{-i}d$).

For the high frequency we should do something similar: one would like to localize $|\eta|$ in intervals of size $2^{-j}d$. The only problem we encounter is that if $i << j$ and $d$ small we may see the curvature of the circle and then the support of $g_{l^{\perp},n}$ cannot be approximated by a rectangle.   

In order to fix this we chose $g_{l^{\perp},n}$ to be the part of $g_{l}$ localized in  $A^{j}_{l^{\perp},n}=\{\eta \in A_{l^{\perp}}: \eta \cdot v_{l^{\perp}} \in [2^{-j}d(n-\q),2^{-j}d(n+\q)] \}$; we denoted by $v_{l^{\perp}}=(\cos(l^{\perp}2^{-i-j}d \frac{\pi}{2}),\sin(l^{\perp}2^{-i-j}d \frac{\pi}{2}))$ and $v_{l}=(\cos(l2^{-i-j}d \frac{\pi}{2}),\sin(l2^{-i-j}d \frac{\pi}{2}))$ (we will need the second one later).

This way the supports of $A^{i}_{l,m}$ and $B^{j}_{l^{\perp},n}$ are rectangles of sizes $2^{-i}d \times 2^{-j}d$.

The crucial property is that the sum sets of the supports, namely $A^{i}_{l,m}+A^{j}_{n,l^{\perp}}=\{\xi+\eta:\xi \in A^{i}_{l,m} \ \mbox{and} \ \eta \in  A^{j}_{n,l^{\perp}}\}$ are disjoint with respect to the pair $(m,n)$. This is mainly because the sum set $A^{i}_{l,m}+A^{j}_{n,l^{\perp}}$ is approximately a rectangle of sizes $2^{-i}d \times 2^{-j}d$ and whose center has coordinates $(2^{-j}dn,2^{-i}dm)$ with respect to the base $(v_{l},v_{l^{\perp}})$. Let's denote by $h_{m,n}$ the part of $h$ which is supported in this set (more precisely the projection of the support on the $\xi$ space should be supported there). Hence we can write:

$$(f \delta_{\tau=\xi^{2}} * g \delta_{\tau=\xi^{2}})h=$$

$$\sum_{l \in J_{i,j,d}} \sum_{m} \sum_{n} \int f_{l,m}(\xi) g_{l^{\perp},n}(\eta) h_{m,n}(\xi+\eta, \xi^{2}+\eta^{2}) \sqrt{1+4\xi^{2}}  \sqrt{1+4\eta^{2}} d\xi d\eta $$

\noindent
and then, for fixed $l$, we can estimate:

$$|(f_{l} \delta_{\tau=\xi^{2}} * g_{l^{\perp}} \delta_{\tau=\xi^{2}})h| \les  \left( \sum_{m} ||f_{l,m}||^{2}_{L^{2}(P)} \right)^{\q} \left( \sum_{n} ||g_{l^{\perp},n}||^{2}_{L^{2}(P)} \right)^{\q} \cdot$$

$$\left( \sum_{m,n} \int h^{2}_{m,n}(\xi+\eta,\xi^{2}+\eta^{2}) \sqrt{1+4\xi^{2}}  \sqrt{1+4\eta^{2}} d\xi d\eta \right)^{\q}$$

In order to estimate the remaining integral, we introduce the change of variables $\eta \to (\r,\t)$ followed by $(\xi_{1},\xi_{2},\r) \to (\z_{1},\z_{2},\z_{3})$:

\begin{eqnarray} \label{sys}
  \left\{
        \begin{array}{rl}
		 \xi_{1}+\r\cos{\t}=\z_{1}  \\
		 \xi_{2}+\r\sin{\t}=\z_{2}  \\
		 \xi_{1}^{2}+\xi^{2}_{2}+\r^{2}=\z_{3}
	\end{array}\right.
\end{eqnarray}

The Jacobian for the first transformation is $d \eta_{1} d \eta_{2}=\r d\r d\t$ and for the second $d\z_{1} d\z_{2} d\z_{3}=(\r-\xi_{1}\cos{\t}-\xi_{2}\sin{\t})d\xi_{1} d\xi_{2} d\r=\frac{\r^{2}-\xi \cdot \eta}{\r}d\xi_{1} d\xi_{2} d\r \approx \r d\xi_{1} d\xi_{2} d\r$; here we have used the fact that $|\xi \cdot \eta| \leq d \leq \frac{\r^{2}}{4}$. Therefore the above integral becomes:

$$\int h^{2}_{m,n}(\xi+\eta,\xi^{2}+\eta^{2}) \sqrt{1+4\xi^{2}}  \sqrt{1+4\eta^{2}} d\xi d\eta \approx$$

$$2^{i+j}\int h^{2}_{m,n}(\z_{1},\z_{2},\z_{3}) d\z_{1} d\z_{2} d \z_{3} d\t  \leq d ||h_{m,n}||^{2}$$

In the last estimate we have used the fact that we integrate over a domain where $\Delta \t \approx 2^{-i-j}d$. Hence we can conclude the above computation with:

$$|(f_{l} \delta_{\tau=\xi^{2}} * g_{l^{\perp}} \delta_{\tau=\xi^{2}})h| \les  d^{\q}||f_{l}||_{L^{2}(P)}  ||g_{l^{\perp}}||_{L^{2}(P)} ||h||_{L^{2}}$$

In the end we perform the summation with respect to $l$ to obtain:

$$|(f \delta_{\tau=\xi^{2}} * g \delta_{\tau=\xi^{2}})h| \les \sum_{l} d^{\q} ||f_{l}||_{L^{2}(P)}  ||g_{l^{\perp}}||_{L^{2}(P)} ||h||_{L^{2}} \les$$

$$ d^{\q} \left( \sum_{l} ||f_{l}||^{2}_{L^{2}(P)} \right)^{\q} \left( \sum_{l} ||g_{l^{\perp}}||^{2}_{L^{2}(P)} \right)^{\q} ||h||_{L^{2}} \les d^{\q} ||f||_{L^{2}(P)} ||g||_{L^{2}(P)} ||h||_{L^{2}}$$

Since this holds true for any $h \in L^{2}$ supported in $|\tau-\xi^{2}| \leq d$, we can conclude with the result of the Proposition.

{\mathversion{bold} $f \delta_{\bar{P}} * g \delta_{P}$}

We start the argument in a similar way. We test the convolution against a $h \in L^{2}$ supported in $|\tau-\xi^{2}| \leq d$:

$$(f \delta_{\tau=-\xi^{2}} * g \delta_{\tau=\xi^{2}})h=$$

$$\int f(\xi) g(\eta) h(\xi+\eta, -\xi^{2}+\eta^{2}) \sqrt{1+4\xi^{2}}  \sqrt{1+4\eta^{2}} d\xi d\eta $$

Since $h$ is supported in a region $|\tau-\xi^{2}| \leq d$ we need the following condition on the variables inside the integral: $|(\xi+\eta)^{2}-(-\xi^{2}+\eta^{2})| \leq d$ or $2|\xi||\eta+\xi|\cos{\t} \leq d$ where $\t$ is the angle between $\xi$ and $\eta+\xi$. Hence $\cos{\t} \leq 2^{-i-j}d$ which implies $|\t-\frac{\pi}{2}| \leq 2^{-i-j}d$. This suggests decomposing:

$$f=\sum_{l \in J_{i,j,d}} f_{l}  \ \ \ \mbox{and} \ \ \ h=\sum_{l \in J_{i,j,d}} h_{l}$$

\noindent
such that $f_{l}$ is as before and $h_{l}$ is the part of $h$ whose support, when projected on the $\xi$ plane, is included in $A_{l}$. Then the argument continues as before with $h$ taking the place of $g$ and vice-versa. 

{\mathversion{bold} $f \delta_{P} * g \delta_{\bar{P}}$}

If $i \leq j-2$ then the convolution is localized in region with $\tau \leq 0$, hence outside the region with $|\tau-\xi^{2}| \leq d$. If $j-1 \leq i \leq j$ then this is similar to the case {\mathversion{bold} $f \delta_{\bar{P}} * g \delta_{P}$}.

\end{proof}

\vspace{.2in}

\subsection{Bilinear estimates on dyadic regions}
\noindent
\vspace{.1in}

 For a bilinear estimate we use the notation:

$$ \mathcal{X} \cdot \mathcal{Y} \rightarrow \mathcal{Z}$$

\noindent
which means that we seek for an estimate $||B(u,v)||_{\mathcal{Z}} \leq C ||u||_{\mathcal{X}} \cdot ||v||_{\mathcal{Y}}$. Here the constant $C$ may depend on some variables, like the frequency where the functions are localized.

A standard way of writing down each case looks like:

\vspace{.1in}

{\mathversion{bold} $ X_{i, d_{1}}^{0,\q} \cdot X_{j, d_{2}}^{0,\q} \rightarrow  X_{j, d_{3}}^{0,-\q}$  }

\vspace{.1in}

This means that for $u \in X_{i, d_{1}}^{0,\q}$ and $v \in X_{j, d_{2}}^{0,\q}$ we estimate the part of $B(u,v)$ (or $\tilde{B}(u,v)$) whose Fourier transform is supported in $A_{j, d_{3}}$. Formally we estimate $\mathcal{F}^{-1}(\chi_{A_{j,d_{3}}}\mathcal{F}(B(u,v)))$. This is going to be the only kind of ``abuse'' in notation which we make throughout the paper, i.e. considering $||B(u,v)||_{X_{j,d_{3}}^{s,\q}}$ even if $\mathcal{F}(B(u,v))$ is not supported in $A_{j,d_{3}}$. We choose to do this so that we do not have to relocalize every time in $A_{j,d_{3}}$.

 Sometimes we prove estimates via duality or conjugation: 

$$X \cdot Y \rightarrow Z \Longleftrightarrow X \cdot (Z)^{*} \rightarrow (Y)^{*}  \ \ \ \mbox{and} \ \ \   X \cdot Y \rightarrow Z \Longleftrightarrow \bar{X} \cdot \bar{Y} \rightarrow \bar{Z}$$

\begin{p3} \label{u1}

Assume $1 \leq i \leq j$. Then we have the estimates:

\beq \label{b1}
||B(u,v)||_{X^{0,-\q}_{j,d_{3}}} \les 2^{\frac{3i}{2}} (\max{(d_{2},d_{3})})^{-\q} ||u||_{X^{0,\q}_{i,d_{1}}} ||v||_{X^{0,\q}_{j,d_{2}}}
\eeq

\beq \label{b2}
||B(u,v)||_{X^{0,-\q}_{j,d_{3}}} \les 2^{\frac{i+j}{2}} ||u||_{X^{0,\q}_{i,d_{1}}} ||v||_{X^{0,\q}_{j,d_{2}}}
\eeq

\noindent
where the parameters involved are restricted by $i \leq j-5 \Rightarrow d_{1} \leq 2^{i-3}$. 

If $|i-j| \leq 1$ and $k \leq j-1$, then we have the estimates:

\beq \label{b10}
||B(u,v)||_{X^{0,-\q}_{k,d_{3}}} \les 2^{\frac{k}{2}+j} (\max{(d_{1},d_{2})})^{-\q} ||u||_{X^{0,\q}_{i,d_{1}}} ||v||_{X^{0,\q}_{j,d_{2}}}
\eeq

\beq \label{b11}
||B(u,v)||_{X^{0,-\q}_{k,d_{3}}} \les 2^{\frac{3j-k}{2}} ||u||_{X^{0,\q}_{i,d_{1}}} ||v||_{X^{0,\q}_{j,d_{2}}}
\eeq

\noindent
where the parameters are restricted by $k \leq j-5 \Rightarrow d_{3} \leq 2^{k-2}$.

All of the above estimates hold true, with the same restrictions, if $B(u,v)$ is replaced by $B(\bar{u},v)$ or $B(u,\bar{v})$. 

\end{p3}

\begin{proof}

We should make some commentaries about the statement above. If $i \leq j-2$, then the result is localized at frequency $\approx 2^{j}$. There is something to estimate only if $k=j, j \pm 1$. The estimates for the case $k=j$ are generic and this is why we choose to list and prove them only.  

It is only when $i=j-1,j$ that we have parts of the result at lower frequencies and then we have to provide estimates for all $k \leq j+1$.

We deal first with the case when we measure the outcome at the high frequency and at the end we deal with the case when we have $i=j-1,j$ and we have to measure the outcome at lower frequencies. 

We need  to transform the estimates on paraboloids in estimates on dyadic pieces. If we localize in a region where $|\xi| \approx 2^{k}$, the parabolas $P_{c}$ make an angle of $\approx 2^{-k}$ with the $\tau$ axis, so we have the following relation between measures:

$$d\xi d\tau \approx 2^{-k} dP_{c} dc$$

If $d \leq 2^{i-3}$ then in $A_{i,d}$ we have $|\xi| \approx 2^{i}$. Therefore for $l \leq i-3$:

\beq \label{e100}
||u||^{2}_{X^{0,\pm \q}_{i,2^{l}}} \approx 2^{\pm(l+i)} \int_{b=2^{l-1}}^{2^{l+1}} ||\hat{u}||^{2}_{L^{2}(P_{(b2^{i})})} db 
\eeq

At this time we are ready to start the estimates.

\vspace{.1 in}

{\mathversion{bold} $ X_{i, d_{1}}^{0,\q} \cdot  X_{j, d_{2}}^{0,\q} \rightarrow  X_{j, d_{3}}^{0,-\q}$}

\vspace{.1 in}

\begin{bfseries} Case 1: $d_{1} \leq 2^{i-3}$ \end{bfseries}

\begin{bfseries} subcase 1.1: $d_{2} \leq 2^{j-3}$ \end{bfseries}

 We can apply the result of (\ref{ge}) to evaluate

$$||\hat{u} * \hat{v}||_{L^{2}} \les \int_{I_{1}} \int_{I_{2}} ||\hat{u} \delta_{P_{b_{1}2^{i}}} * \hat{v} \delta_{P_{b_{2}2^{j}}} ||_{L^{2}} db_{1} db_{2} \les $$

$$ \int_{I_{1}} \int_{I_{2}} 2^{i} ||\hat{u}||_{L^{2}(P_{b_{1}2^{i}})} ||\hat{v}||_{L^{2}(P_{b_{2}2^{j}})} db_{1} db_{2} \les $$

$$2^{i} \left( \int_{I_{1}} (1+b_{1}2^{i})^{-1} db_{1} \right)^{\q} ||u||_{X^{0,\q}_{i,d_{1}}} \left( \int_{I_{2}} (1+b_{2}2^{j})^{-1} db_{2} \right)^{\q} ||v||_{X^{0,\q}_{j,d_{2}}} \approx $$

$$2^{\frac{i-j}{2}} ||u||_{X^{0,\q}_{i,d_{1}}} ||v||_{X^{0,\q}_{j,d_{2}}}$$

Here we used the fact that $I_{1} \approx [\frac{d_{1}}{2},2d_{1}]$ which gives us $\int (1+b_{1}2^{i})^{-1} db_{1} \approx 2^{-i} $. Same thing for the integral with respect to $b_{2}$. (\ref{e100}) gives us:

$$||\tilde{B}(u,v)||_{X^{0,-\q}_{j,d_{3}}} \les (2^{j}d_{3})^{-\q} ||\hat{u} * \hat{v}||_{L^{2}} \les 2^{\frac{i}{2}-j} d_{3}^{-\q} ||u||_{X^{0,\q}_{i,d_{1}}} ||v||_{X^{0,\q}_{j,d_{2}}}$$ 

If we use the same approach, but use (\ref{ge1}) instead, we obtain the estimate (\ref{b2}). The use of (\ref{ge1}) requires some restrictions on the range of parameters involved. In our case this translates into $d_2,d_3 \leq 2^{i-2}$. This is fine, since if $\max{(d_2,d_3)} \geq 2^{i-2}$, then the estimate (\ref{b1}) is stronger, hence we do not have to deal with additional cases.

\begin{bfseries} subcase 1.2: $d_{3} \leq 2^{j-3}$ \end{bfseries}

This estimate for this case can be deduced by duality from the estimate:

$$ X_{i, d_{1}}^{0,\q} \cdot  \bar{X}_{j, d_{3}}^{0,\q} \rightarrow  \bar{X}_{j, d_{2}}^{0,-\q} \Leftrightarrow  \bar{X}_{i, d_{1}}^{0,\q} \cdot  X_{j, d_{3}}^{0,\q} \rightarrow  X_{j, d_{2}}^{0,-\q} $$

The proof of the last estimate is treated in Subcase 1.1 bellow.

\begin{bfseries} subcase 1.3: $d_{2},d_{3} \geq 2^{j-2}$  \end{bfseries}

In this case we use a much simpler argument. For reference we call this the $L^{1} * L^{2} \rightarrow L^{2}$ argument. It goes as follows:

$$||\hat{u}||_{L^{1}} \les 2^{\frac{3}{2}i} d_{1}^{\q} ||\hat{u}||_{L^{2}} \approx 2^{i} ||u||_{X_{i,d_{1}}^{0,\q}} $$

 Then we continue with:

$$||\tilde{B}(u,v)||_{X^{0,-\q}_{j,d_{3}}} \approx  2^{-j} ||\hat{u} * \hat{v} ||_{L^{2}(A_{j,d_{3}})} \les $$

$$ 2^{-j} ||\hat{u}||_{L^{1}} \cdot ||\hat{v}||_{L^{2}} \les 2^{i-2j} ||u||_{X^{0,\q}_{j,d_{1}}} ||v||_{X^{0,\q}_{j,d_{2}}}$$

Notice that, since $i \leq j$, the first estimate is stronger then the second, hence we do not have anything else to prove.

\begin{bfseries} Case 2: $d_{1} \geq 2^{i-2}$  \end{bfseries}

 We have to deal only with the case $i \geq j-5$. We get this estimate by duality from:

$$\bar{X}_{j, d_{3}}^{0,\q} \cdot  X_{j, d_{2}}^{0,\q} \rightarrow  \bar{X}_{i, d_{1}}^{0,-\q} \Leftrightarrow  X_{j, d_{3}}^{0,\q} \cdot  \bar{X}_{j, d_{2}}^{0,\q} \rightarrow  X_{i, d_{1}}^{0,-\q}$$

The last estimate is treated in the next group of estimates. 

\vspace{.1 in}

{\mathversion{bold} $ \bar{X}_{i, d_{1}}^{0,\q} \cdot X_{j, d_{2}}^{0,\q} \rightarrow X_{j, d_{3}}^{0,-\q}$}

\vspace{.1 in}

\begin{bfseries} Case 1: $d_{1} \leq 2^{i-3}$ \end{bfseries}

\begin{bfseries} subcase 1.1: $d_{2} \leq 2^{j-3}$ \end{bfseries}

This case is totally similar to Subcase 1.1 in the first estimate because we have all the necessary ingredients.

\begin{bfseries} subcase 1.2: $d_{3} \leq 2^{j-3}$ \end{bfseries}

This estimate for this case can be deduced by duality from the estimate:

$$ \bar{X}_{i, d_{1}}^{0,\q} \cdot  \bar{X}_{j, d_{3}}^{0,\q} \rightarrow  \bar{X}_{j, d_{2}}^{0,-\q} \Leftrightarrow  X_{i, d_{1}}^{0,\q} \cdot  X_{j, d_{3}}^{0,\q} \rightarrow  X_{j, d_{2}}^{0,-\q}$$

Last estimate was proved by itself in Subcase 1.1 in the first group of estimates.

\begin{bfseries} subcase 1.3: $d_{2},d_{3} \geq 2^{j-2}$  \end{bfseries}

Use the $L^{1} * L^{2} \rightarrow L^{2}$ argument.

\begin{bfseries} Case 2: $d_{1} \geq 2^{i-2}$ \end{bfseries}

 We work in the hypothesis $i \geq j-5$. By duality we get the estimate from

$$ \bar{X}_{j, d_{3}}^{0,\q} \cdot \bar{X}_{j, d_{2}}^{0,\q} \rightarrow  X_{i, d_{1}}^{0,-\q}$$

This last estimate can be easily treated as if it were $ \bar{X}_{j, d_{3}}^{0,\q} \cdot \bar{X}_{j, d_{2}}^{0,\q} \rightarrow  \bar{X}_{i, d_{1}}^{0,-\q}$, since $d_{1} \geq 2^{i-2}$. The conjugate of this estimate has been treated before.  

\vspace{.1 in}

{\mathversion{bold} $ X_{i, d_{1}}^{0,\q} \cdot   \bar{X}_{j, d_{2}}^{0,\q} \rightarrow  X_{j, d_{3}}^{0,-\q}$}

\vspace{.1 in}

If $j-5 \leq i \leq j$ then the estimate is similar to the one in {\mathversion{bold} $ \bar{X}_{i, d_{1}}^{0,\q} \cdot  X_{j, d_{2}}^{0,\q} \rightarrow  X_{j, d_{3}}^{0,-\q}$}. If $i \leq j-5$, then we have the following cases:

\begin{bfseries} Case 1: $d_{2},d_{3} \leq 2^{j-3}$ \end{bfseries}

This is incompatible since functions in $\bar{X}_{i, d_{1}}^{0,\q}$ have their Fourier transform supported in a region with $\tau < 0$ and functions in $X_{j, d_{3}}^{0,-\q}$ have their Fourier transform supported in a region with $\tau > 0$; an easy computation shows that, by convolution, the Fourier transform of a function in $X_{i, d_{1}}^{0,\q}$ cannot move the first support to the second one. 

\begin{bfseries} Case 2: $\min{(d_{2},d_{3})} \geq 2^{j-3}$ \end{bfseries}

\begin{bfseries} subcase 2.1: $d_{3} \geq 2^{2j-3} \geq d_{2}$ \end{bfseries}

This case can be treated the same way like {\mathversion{bold} $ X_{i, d_{1}}^{0,\q} \cdot  X_{j, d_{2}}^{0,\q} \rightarrow  X_{j, d_{3}}^{0,-\q}$}, just that we use this time the estimate for $f \delta_{P_{c_{1}}} * g \delta_{\bar{P}_{c_{2}}}$. Notice that the condition $d_{3} \geq 2^{j-3}$ implies that we have to deal only with the estimate (\ref{b1}), since it becomes stronger than the estimate (\ref{b2}).

\begin{bfseries} subcase 2.2: $d_{2} \geq 2^{2j-3} \geq d_{3}$ \end{bfseries}

This case can be obtained by duality from  {\mathversion{bold} $ X_{i, d_{1}}^{0,\q} \cdot  \bar{X}_{j, d_{3}}^{0,\q} \rightarrow  X_{j, d_{2}}^{0,-\q}$}, which is similar to the above estimate.

\begin{bfseries} subcase 2.3: $\max{(d_{2},d_3)} \geq 2^{2j-3}$ \end{bfseries}

Use the $L^{1} * L^{2} \rightarrow L^{2}$ argument.

\vspace{.1 in}

\begin{bfseries} High - High interactions with output at low frequencies  \end{bfseries}

\vspace{.1 in}

We have to deal with estimates of type $ X_{i, d_{1}}^{0,\q} \cdot  X_{j, d_{2}}^{0,\q} \rightarrow   X_{k, d_{3}}^{0,-\q}$, for $i=j-1,j$ and $k \leq j+1$. The estimate for the case $i=j$ is generic, hence we will work only this one out. In order to see easier the duality, we choose to replace $k$ by $i$ (this new $"i"$ is different the the one before) and look for an estimate of type

{\mathversion{bold} $  X_{j, d_{1}}^{0,\q} \cdot  X_{j, d_{2}}^{0,\q} \rightarrow   X_{i, d_{3}}^{0,-\q}$}

\vspace{.1 in}

Conjugation and duality give us:

$$ X_{i, d_{3}}^{0,\q} \cdot \bar{X}_{j, d_{2}}^{0,\q} \rightarrow  X_{j, d_{1}}^{0,-\q} \Rightarrow  \bar{X}_{i, d_{3}}^{0,\q} \cdot  X_{j, d_{2}}^{0,\q} \rightarrow  \bar{X}_{j, d_{1}}^{0,-\q} \Rightarrow X_{j, d_{1}}^{0,\q} \cdot X_{j, d_{2}}^{0,\q} \rightarrow  X_{i, d_{3}}^{0,-\q}$$

\noindent
and this is enough to justify the estimate; with one exception though: $i \geq j-5$ and $d_{3} \geq 2^{i-2}$. This exception is treated in the next two cases.

\begin{bfseries} Case 1: $d_{1}, d_{2} \leq 2^{j-2}$ \end{bfseries}

The argument is similar to Subcases 1.1 in the previous estimates. Making use of (\ref{ge}) we get:

$$||\hat{u} * \hat{v}||_{L^{2}} \les \int_{I_{1}} \int_{I_{2}} ||\hat{u} \delta_{P_{b_{1}2^{i}}} * \hat{v} \delta_{P_{b_{2}2^{j}}} ||_{L^{2}} db_{1} db_{2} \les$$

$$ \int_{I_{1}} \int_{I_{2}} 2^{j} ||\hat{u}||_{L^{2}(P_{b_{1}2^{j}})} ||\hat{v}||_{L^{2}(P_{b_{2}2^{j}})} db_{1} db_{2} \les  ||u||_{X^{0,\q}_{j,d_{1}}} ||v||_{X^{0,\q}_{j,d_{2}}}$$

Next

$$||\tilde{B}(u,v)||_{X^{0,-\q}_{i,d_{3}}} \approx  (2^{2j})^{-\q} ||\tilde{B}(u,v)||_{L^{2}(A_{i,d_{3}})} \les 2^{-j} ||u||_{X^{0,\q}_{j,d_{1}}} ||v||_{X^{0,\q}_{j,d_{2}}}$$

\noindent
where we use the fact that $d_{3} \geq 2^{i-2} \geq 2^{j-7}$. 

\begin{bfseries} Case 2: $\max{(d_{1}, d_{2})} \geq 2^{j-2}$ \end{bfseries}

This case is similar to Subcases 1.3 in the estimate {\mathversion{bold} $ X_{i, d_{1}}^{0,\q} \cdot  X_{j, d_{2}}^{0,\q} \rightarrow  X_{j, d_{3}}^{0,-\q}$} and uses the trivial $L^{1} * L^{2} \rightarrow L^{2}$ argument. 

\vspace{.1 in}

{\mathversion{bold} $X_{j, d_{1}}^{0,\q} \cdot  \bar{X}_{j, d_{2}}^{0,\q} \rightarrow   X_{i, d_{3}}^{0,-\q}$}

\vspace{.1 in}

In the same way as above, duality gives us the estimates as claimed in the Theorem. With the same exception: $i \geq j-5$ and $d_{3} \geq 2^{i-2}$ ! The exception is treated we proceed as in the Case 1 and Case 2 above. 

\end{proof}

\vspace{.1in}

\subsection{Bilinear estimates in $X^{s,\q}$ involving decay}
\noindent
\vspace{.1in}

The results in Proposition \ref{u1} indicate that we can not recover a full derivative in the case when, see (\ref{b1}) and (\ref{b2}), $d_2$ and $d_3$ are small. We also do not have a complete range of bilinear estimates, see the restriction $d_1 \leq 2^{i-3}$ there. This is an indication that we have to bring the additional decay structure in order to complete the section.

\begin{p3} \label{u2}

a) Assume $i \leq j-5$, $d_{1} \leq 2^{i-2}$ and $d_{2},d_{3} \leq 2^{i-5}$. Then we have the following estimates:

\beq \label{bdec1}
||B(u,v)||_{X^{0,-\q}_{j,d_{3}}} \les 2^{2i} ||u||_{\mathcal{D}X^{0,\q}_{i,d_{1}}} ||v||_{X^{0,\q}_{j,d_{2}}}
\eeq

b) Assume $i \leq j-5$ and $d_{1} \geq 2^{i-2}$. Then we have the following estimates:

\beq \label{bdec2}
||B(u,v)||_{X^{0,-\q}_{j,d_{3}}} \les 2^{2i} ||u||_{\mathcal{D}X^{0,\q}_{i,d_{1}}} ||v||_{X^{0,\q}_{j,d_{2}}}
\eeq

The above estimates holds true if $B(u,v)$ is replaced by $B(\bar{u},v)$ or $B(u,\bar{v})$. 

c) Assume $|i-j| \leq 1$, $k \leq j-5$ and $d_{3} \geq 2^{i-2}$. Then we have the following estimates:

\beq \label{bdec3}
||B(u,v)||_{X^{0,-\q}_{k,d_{3}}} \les 2^{k+j} 2^{\frac{j-k}{2}}||u||_{\mathcal{D}X^{0,\q}_{i,d_{1}}} ||v||_{X^{0,\q}_{j,d_{2}}}
\eeq

The above estimates holds true if $B(u,v)$ is replaced by $B(\bar{u},v)$ or $B(u,\bar{v})$.

\end{p3}

\begin{proof}

a) It is enough to estimate $\langle f \cdot g, h \rangle_{L^{2}}=\langle \hat{f} * \hat{g}, \hat{h} \rangle_{L^{2}}$ for $f \in \mathcal{D}X_{i,d_{1}}^{0,\q}$, $g \in X_{j,d_{2}}^{0,\q}$ and $h \in (X_{j,d_{3}}^{0,-\q})^{*}=\bar{X}^{0,\q}_{j,d_{3}}$. 

We define $\Theta_{p}=\{\eta \in \Xi, |\eta| \approx 2^{j}: \arg{\eta} \in [(p-\q)2^{i-j}\frac{\pi}{2},(p+\q)2^{i-j}\frac{\pi}{2}]\}$ and by $g_{p}$ we denote the part of $f$ whose Fourier transform is localized in $\Theta_{p}$. The size of $\hat{g}_{p}$ in the angular direction is $\approx 2^{i}$, therefore the interactions $\hat{f} * \hat{g}_{p}$ is essentially localized in  $\Theta_{p}$. Therefore we have:

$$\langle f \cdot g, h \rangle_{L^{2}}=\sum_{p} \langle f \cdot g_{p}, h_{p} \rangle_{L^{2}} $$

For each $p$, we define orthonormal basis $x_{p}^1,x_{p}^2$ in $\R^{2}$ by $x_{p}^{1}=(\cos{(p2^{i-j}\frac{\pi}{2})},\sin{(p2^{i-j}\frac{\pi}{2})})$ and $x_{p}^{2}=(-\sin{(p2^{i-j}\frac{\pi}{2})},\cos{(p2^{i-j}\frac{\pi}{2})})$. We denote by $\xi_{p}^{1}, \xi_{p}^2$ the corresponding basis on the Fourier side.

For each $m \in Z$ we define the rectangles $R^{m}_{p} \subset \R^{2}$ centered at $(m_{1}2^{j},m_{2}2^{i})$ with respect to the basis $(x^{p},y^{p})$ and of sizes $2^{j} \times 2^{i}$ ($x^{p} \times y^{p}$ directions).

$\hat{g}_{p}$ is supported in a "curved" parallelepiped whose sizes are larger than the dual sizes of $Q^{m}_{p} \times \R$; the size of its support in the direction of $\xi_p^1$ is $d_{2}$, hence we can conclude:

$$\sum_{m} ||g_{p}||_{L^{2}_{y_{p},t}L^{\infty}_{x^{p}}(R_{p}^{m} \times \R)}^{2} \les d_{2} ||g_{p}||^{2}_{L^{2}}$$

Similarly for $h_{p,l}$ (here it is key that it has the same angular localization):

$$\sum_{m} ||\bar{h}_{p}||_{L^{2}_{y_{p},t}L^{\infty}_{x^{p}}(R_{p}^{m} \times \R)}^{2} \les d_{3} ||\bar{h}_{p}||^{2}_{L^{2}}$$

This gives us an $l^{1}_{m} (L^{1}_{y^{p},t}L^{\infty}_{x^{p}}(R^{m}_{p} \times \R))$ structure for $g_{p}\bar{h}_{p}$. Therefore we need to estimate $f$ in $l^{\infty}_{m} (L^{\infty}_{y^{p},t}L^{1}_{x^{p}}(R^{m}_{p} \times \R))$.

Each $R_{p}^{m}$ can be written as $R_{p}^{m} \subset \cup_{\bar{m}:Q_{i}^{\bar{m}} \cap L^{k}_{i} \ne \emptyset} Q_{i}^{\bar{m}}$. We can be more precise: the line $L$ generating $L^{k}_{i}$ goes in the direction of $x^{p}$ and we can restrict indexes $\bar{m}$ in a set of cardinality $\approx 2^{j-i}$. Then we have:

$$||f||_{L^{\infty}_{y^{p},t}L^{1}_{x^{p}}(R^{m}_{p} \times \R)} \les \sum_{\bar{m}} ||f||_{L^{\infty}_{y^{p},t}L^{1}_{x^{p}}(Q^{\bar{m}}_{i} \times \R)} \les $$

$$\sum_{\bar{m}} \sum_{\xi} \sum_{k} ||f_{\xi,\xi^{2}+k}||_{L^{\infty}_{y^{p},t}L^{1}_{x^{p}}(Q^{\bar{m}}_{i} \times \R)} \les  2^{i} \sum_{\bar{m}} \sum_{\xi} \sum_{k} ||f_{\xi,\xi^{2}+k}||_{L^{\infty}(Q^{\bar{m}}_{i} \times \R)} \les$$

$$ 2^{i} (2^{3i}d_{1})^{\q} \sum_{\bar{m}} \left(\sum_{\xi} \sum_{k} ||f_{\xi,\xi^{2}+k}||^{2}_{L^{\infty}(Q^{\bar{m}}_{i} \times \R)} \right)^{\q}$$

Now we use (\ref{dd2}) to obtain:

$$||f||_{L^{\infty}_{y^{p},t}L^{1}_{x^{p}}(R^{m}_{p} \times \R)} \les (2^{3i}d_{1})^{\q} ||f||_{\mathcal{D}_{i}L^{2}} \les 2^{i} ||f||_{\mathcal{D}X^{0,\q}_{i,d_{1}}}$$

Hence we obtain the estimate:

$$|\langle f \cdot g_{p}, h_{p} \rangle_{L^{2}}| \les 2^{i} ||f||_{\mathcal{D}X^{0,\q}_{i,d_{1}}} d_{2}^{\q} ||g_{p}||_{L^{2}} d_{3}^{\q} ||\bar{h}_{p}||_{L^{2}}$$

Summing wit respect to $p$ we obtain:

$$|\langle f \cdot g, h \rangle_{L^{2}}| \les 2^{i} ||f||_{\mathcal{D}X^{0,\q}_{i,d_{1}}} d_{2}^{\q} ||g||_{L^{2}} d_{3}^{\q} ||\bar{h}||_{L^{2}} \approx 2^{i-j}||f||_{\mathcal{D}X^{0,\q}_{i,d_{1}}} ||g||_{X^{0,\q}_{j,d_{2}}} ||h||_{\bar{X}^{0,\q}_{j,d_{3}}}$$

This translates into:

$$||\tilde{B}(f,g)||_{X_{j,d_{3}}^{0,\q}} \les 2^{i-j} ||f||_{\mathcal{D}X^{0,\q}_{i,d_{1}}} ||g||_{X^{0,\q}_{j,d_{2}}}$$

The principle in (\ref{bb}) gives us (\ref{bdec1}).

One can easily see from the above argument that we can easily carry on the same proof for the case when we deal with $B(\bar{f},g)$. As about $B(f,\bar{g})$, this was written just for the sake of completeness, since there is nothing to estimate there. The supports of high frequencies $\bar{A}_{j,d_{2}}$ and $A_{j,d_{3}}$ are too far away to be linked via the convolution with the low frequency since we are in the case $i \leq j-5$.

b) The proof for the case $d_{1} \geq 2^{i-2}$, $d_{2},d_{3} \leq 2^{j-2}$ is completely similar of the previous one. We work with $g_{\xi,l}$ instead of $g_{\xi,\xi^{2}+k}$ and use (\ref{dd6}) instead of (\ref{dd2}).

In the case $d_{2} \leq 2^{j-2}$ and $d_{3} \geq 2^{j-2}$ we estimate directly since it turns out that no decay is needed. We use the $L^{1} * L^{2} \to L^{2}$ argument:

$$||f \cdot g||_{X^{0,\q}_{j,d_{3}}} \approx 2^{-j}||\hat{f}*\hat{g}||_{L^{2}} \les 2^{-j}||\hat{f}||_{L^{1}} ||\hat{g}||_{L^{2}} \les $$

$$2^{2i-j} ||\hat{f}||_{L^{2}}  ||g||_{L^{2}} \les 2^{i-j} ||f||_{X^{0,\q}_{i,d_{1}}} ||g||_{X_{j,d_{2}}^{0,\q}}$$

This is the estimate for $\tilde{B}$ which implies the one for $B$, see (\ref{bb}).

In the case when $d_{2} \geq 2^{j-2}$ but $d_{3} \leq 2^{j-2}$ we can obtain the result via duality from $\mathcal{D}X^{0,\q}_{i,d_{1}} \cdot \bar{X}^{0,\q}_{j,d_{3}} \to \bar{X}^{0,\q}_{j,d_{2}}$. This last estimate can be obtained via a similar argument as the one above.

If $d_{2},d_{3} \geq 2^{j-2}$ then the estimate can be obtained via the $L^{1} * L^{2} \to L^{2}$ argument.

As we remark at the end of the proof of a), we can obtain easily the estimates for $B(\bar{u},v)$ and $B(u,\bar{v})$.

c) We will get this result from the one in (\ref{bdec1}). It is enough to estimate $\langle f \cdot g, h \rangle_{L^{2}}$ for $f \in X_{i,d_{1}}^{0,\q}$, $g \in X_{j,d_{2}}^{0,\q}$ and $h \in (X_{k,d_{3}}^{0,-\q})^{*}=\bar{X}^{0,\q}_{k,d_{3}}$. Equivalently, we can estimate $\langle \bar{h} \cdot g, \bar{f} \rangle_{L^{2}}$ for $\bar{h} \in X_{k,d_{3}}^{0,\q}$, $g \in X_{j,d_{2}}^{0,\q}$ and $\bar{f} \in X^{0,\q}_{i,d_{1}}$. This can be treated the same way as we did with (\ref{bdec2}), just that we do not have the $\mathcal{D}$ structure for $f$ anymore.  

 A careful look at the complete argument for (\ref{bdec1}) shows that the lack of this structure requires a factor of $2^{\frac{j-k}{2}}$ to be added to the estimates. Other than that the argument will be carried the same way as the one in part a), just that at the very end there is an additional correction which should be made: the factor of $2^{2k}$ which would have been obtained in part a) is replaced by $2^{k+j}$ since we have gradients on the two high frequencies and not one gradient on the low and one on the high frequency as in part a). With these observation, we completed the argument for (\ref{bdec3}).

\end{proof}

We can put the results of Proposition \ref{u1} and \ref{u2} together to obtain:

\begin{p3} \label{u3}

a)If $i \leq j$, we have the following estimates:

\beq \label{bee1}
||B(u,v)||_{X^{0,-\q}_{j,d_{3}}} \les 2^{2i} ||u||_{\mathcal{D}_{i}X^{0,\q}_{i,d_{1}}} ||v||_{X^{0,\q}_{j,d_{2}}}
\eeq

The above estimates holds true if $B(u,v)$ is replaced by $B(\bar{u},v)$ or $B(u,\bar{v})$. 

b) Assume $|i-j| \leq 1$. Then we have the following estimates:

\beq \label{bee2}
||B(u,v)||_{X^{0,-\q}_{k,d_{3}}} \les 2^{\frac{k+3j}{2}} ||u||_{\mathcal{D}_{i}X^{0,\q}_{i,d_{1}}} ||v||_{X^{0,\q}_{j,d_{2}}}
\eeq

The above estimates holds true if $B(u,v)$ is replaced by $B(\bar{u},v)$ or $B(u,\bar{v})$.

\end{p3}

\vspace{.1in}

\subsection{Abstract result} \label{ar}
\noindent
\vspace{.1in}

Before we turn to completing the bilinear estimates by proving the conservation of decay, we need to prepare some theoretical facts. 

We fix $\l$ positive. The arguments bellow are independent of the size of $\l$, hence, later on, we have to freedom to apply the results we obtain for various values of $\l$.  

Let $(Q^{m})_{m \in \Z^{n}}$ be a partition of $\R^{n}$ (physical space) in disjoint cubes of size $\l$. We assume that $Q^{m}$ is centered at $ \l m$. Similarly, let  $(R^{p})_{p \in \Z^{n}}$ be  a partition of $\R^{n}$ (frequency space) in disjoint cubes of size $\l^{-1}$. We assume that $R^{p}$ is centered at $ \l^{-1} p$. 

Let $\chi_{Q^{m}}$ be a smooth approximation of the characteristic function of $Q^{m}$ in the following sense: $\chi_{Q^{m}}$ is a non-negative function essentially supported in $Q^{m}$ such that $||\frac{\partial^{\a}\chi_{Q^{m}}}{\partial x^{\a}}||_{L^{\infty}} \les \l^{-|\a|}, \forall \a \in \N^{n}$. We also want:

$$\sum_{m} \chi_{Q^{m}}(x)=1 \ \ \forall x \in R^{n}$$

We say that $\tilde{\chi}_{Q^{m}}$ is a generalized characteristic function of $Q^{m}$ if $\tilde{\chi}_{Q^{m}}$ is essentially supported in $ Q^{m}$ and $||\frac{\partial^{\a}\tilde{\chi}_{Q^{m}}}{\partial x^{\a}}||_{L^{\infty}} \les \l^{-|\a|}, \forall \a \in \N^{n}$. The simplest examples of generalized characteristic functions of $Q^{m}$ are $\l^{|\a|} \frac{\partial^{\a}\tilde{\chi}_{Q^{m}}}{\partial x^{\a}}$ for any $\a \in \N^{n}$. A key property of generalized functions is that:

\beq \label{ps5}
|\tilde{\chi}_{Q^{m}}(x)| \les \sum_{e \in I} \chi_{Q^{m+e}}(x)
\eeq 

\noindent
where $I=\{v \in \Z^{n}: |v_{i}| \leq 1 \forall i=1,..,n \}$.  

In a similar way we define the system $(\chi_{R^{p}})_{p \in \Z^{n}}$; there is only one difference: $||\frac{\partial^{\a}\chi_{R^{p}}}{\partial \xi^{\a}}||_{L^{\infty}} \les \l^{|\a|}, \forall \a \in \N^{n}$. And then we define a generalized characteristic function of $R^{p}$. This time the simplest examples of generalized characteristic functions of $R^{p}$ are $\l^{-|\a|} \frac{\partial^{\a}\tilde{\chi}_{R^{p}}}{\partial \xi^{\a}}$ for any $\a \in \N^{n}$. These generalized characteristic functions enjoy a similar property to (\ref{ps5}).

It is important to make the following convention. While the systems  $(\chi_{Q^{m}})_{m \in \Z^{n}}$ and $(\chi_{R^{p}})_{p \in \Z^{n}}$ are fixed, the generalized characteristic functions can be arbitrary.

\begin{p3} We have the following estimates:

\beq \label{ps1}
||\tilde{\chi}_{Q^{m}}(x) \tilde{\chi}_{R^{p}}(D) \tilde{\chi}_{Q^{m'}}(x) f||_{L^{2}} \les 
\eeq

$$C_{N} \langle m-m'\rangle^{-N} \sum_{|\a|+|\b|= N} ||\tilde{\chi}_{R^{p}}^{\a} (D) \tilde{\chi}_{Q^{m'}}^{\b}(x) f||_{L^{2}}$$

\beq \label{ps4}
\sum_{p} ||\tilde{\chi}_{Q^{m}}(x) \tilde{\chi}_{R^{p}}(D) f||^{2}_{L^{2}} \les C_{N} \sum_{m'} \langle m-m'\rangle^{-N} || \chi_{Q^{m'}}(x) f||^{2}_{L^{2}}
\eeq

\beq \label{ps2}
||\tilde{\chi}_{R^{p}}(D) \tilde{\chi}_{Q^{m}}(x) \tilde{\chi}_{R^{p'}}(D) \chi_{Q^{m'}}(x) f||_{L^{2}} \les 
\eeq

$$C_{N} \langle m-m'\rangle^{-N}  \langle p-p'\rangle^{-N} \sum_{|\a|+|\b| \leq  2N} || \tilde{\chi}_{R^{p}}^{\a} (D) \tilde{\chi}_{Q^{m'}}^{\b}(x) f||_{L^{2}}$$

\noindent
where $\tilde{\chi}_{R^{p}}^{\a}$ is a generalized characteristic function of $R^{p}$ which depends on $\tilde{\chi}_{R^{p}}$ and $\a$ and $\tilde{\chi}_{Q^{m'}}^{\b}$ is a generalized characteristic function of $Q^{m'}$ which depends on $\tilde{\chi}_{Q^{m'}}$ and $\b$.

\end{p3}

Remark. It is not important the exact expression of $\tilde{\chi}_{R^{p}}^{\a}$ or $\tilde{\chi}_{Q^{m'}}^{\b}$ since these terms will be dealt via the estimate (\ref{ps5}). 

\begin{proof}

For (\ref{ps1}) we need to estimate $||\l^{-|\g|} (x-\l m')^{\g} \tilde{\chi}_{R^{p}}(D) \tilde{\chi}_{Q^{m'}}(x) f||_{L^{2}}$. We start with the commutator identity:

$$\l^{-|\g|} (x-\l m')^{\g} \tilde{\chi}_{R^{p}}(D) \tilde{\chi}_{Q^{m'}}(x) f=$$

$$\sum_{\a+\b=\g} \l^{-|\a|} \frac{\partial^{\a} \tilde{\chi}_{R^{p}}} {\partial \xi^{\a}} \l^{-|\b|} (x-\l m')^{\b} \tilde{\chi}_{Q^{m'}}(x) f$$

Then we notice that $\tilde{\chi}^{\b}_{Q^{m'}}=\l^{-|\b|} (x-\l m')^{\b} \tilde{\chi}_{Q^{m'}}$ is a generalized characteristic function of $Q^{m'}$ and $\tilde{\chi}^{\a}_{R^{p}} = \l^{-|\a|} \frac{\partial^{\a} \tilde{\chi}_{R^{p}}}{\partial \xi^{\a}} $ is a generalized characteristic function of $R^{p}$. We denote by $\g^{i}$ the vector in $\N^{n}$ whose $i$'th component is $1$ and the rest are $0$. For $\langle m-m'\rangle \geq 2$ we have:

$$||\langle m-m'\rangle^{N} \tilde{\chi}_{Q^{m}}(x) \tilde{\chi}_{R^{p}}(D) \tilde{\chi}_{Q^{m'}}(x) f||_{L^{2}} \les $$

$$\sum_{i} ||\langle m_{i}-m_{i}'\rangle^{N} \tilde{\chi}_{Q^{m}}(x) \tilde{\chi}_{R^{p}}(D) \tilde{\chi}_{Q^{m'}}(x) f||_{L^{2}} \les$$

$$\sum_{i} ||\l^{-N} (x-\l m')^{N\g^{i}} \tilde{\chi}_{R^{p}}(D) \tilde{\chi}_{Q^{m'}}(x) f||_{L^{2}} \les$$

$$\sum_{i} \sum_{\a+\b=N} ||\tilde{\chi}_{R^{p}}^{\a} (D) \tilde{\chi}_{Q^{m'}}^{\b}(x) f||_{L^{2}}$$

In order to prove (\ref{ps4}) proceed as follows:

$$\sum_{p} ||\tilde{\chi}_{Q^{m}}(x) \tilde{\chi}_{R^{p}}(D) f||^{2}_{L^{2}} \les \sum_{p} \left( \sum_{m'} ||\tilde{\chi}_{Q^{m}}(x) \tilde{\chi}_{R^{p}}(D) \tilde{\chi}_{Q^{m'}}(x) f||_{L^{2}} \right)^{2} \les$$

$$C_{N}^{2} \sum_{p} \left( \sum_{m'} \langle m-m'\rangle^{-N} \sum_{|\a|+|\b|= N} ||\tilde{\chi}_{R^{p}}^{\a} (D) \tilde{\chi}_{Q^{m'}}^{\b}(x) f||_{L^{2}}\right)^{2} \les $$

$$C_{N} \sum_{m'} \langle m-m'\rangle^{-N} \sum_{p} \sum_{|\a|+|\b|= N} ||\tilde{\chi}_{R^{p}}^{\a} (D) \tilde{\chi}_{Q^{m'}}^{\b}(x) f||_{L^{2}}^{2} \les$$

$$C_{N} \sum_{m'} \langle m-m'\rangle^{-N} ||\chi_{Q^{m'}}(x) f||_{L^{2}}^{2} $$

In the last estimate we used twice (once in frequency and once in space) the property (\ref{ps5}).

The proof of (\ref{ps2}) is a direct consequence of (\ref{ps1}) and of its analogue:

$$||\tilde{\chi}_{R^{p}}(x) \tilde{\chi}_{Q^{m}}(D) \tilde{\chi}_{R^{p'}}(x) h||_{L^{2}} \les C_{N} \langle p-p'\rangle^{-N} \sum_{|\a|+|\b|= N} ||\tilde{\chi}_{Q^{m}}^{\a} (D) \tilde{\chi}_{R^{p'}}^{\b}(x) h||_{L^{2}}$$

The proof of this estimate is similar to the one we provided for (\ref{ps1}). If we take in this estimate $h=\tilde{\chi}_{Q^{m'}}f$, we obtain:

$$||\tilde{\chi}_{R^{p}}(x) \tilde{\chi}_{Q^{m}}(D) \tilde{\chi}_{R^{p'}}(x) \tilde{\chi}_{Q^{m'}}f||_{L^{2}} \les $$

$$C_{N} \langle p-p'\rangle^{-N} \sum_{|\a|+|\b|= N} ||\tilde{\chi}_{Q^{m}}^{\a} (D) \tilde{\chi}_{R^{p'}}^{\b}(x) \tilde{\chi}_{Q^{m'}}f||_{L^{2}}$$

Then we apply (\ref{ps1}) for each of the terms $||\tilde{\chi}_{Q^{m}}^{\a} (D) \tilde{\chi}_{R^{p'}}^{\b}(x) \tilde{\chi}_{Q^{m'}}f||_{L^{2}}$ and conclude with the claim in (\ref{ps2}).

\end{proof}

\vspace{.1in}

\subsection{Conservation of decay in bilinear estimates} \label{cd}
\noindent
\vspace{.1in}

We want to warn the reader that this section would be extremely long and tedious if we were to carry out all the computations. This is why will just indicate the main ideas. In principle things should be simple. In the bilinear estimates we used the $\mathcal{D}$ property only on the low frequency, hence the result should inherit the $\mathcal{D}$ property from the high frequency. Which looks reasonable if the interaction is localized at the high frequency too. In the case of high-high to low frequency, there is enough room to transform the  $\mathcal{D}$ structure at high frequency into one at low frequency.

The section is dedicated to proving the following result:

\begin{p3} \label{pdec}

a) If $i \leq j$, we have the following estimates:

\beq \label{be1}
||B(u,v)||_{\mathcal{D}_{j}X^{0,-\q}_{j,d_{3}}} \les 2^{2i} ||u||_{\mathcal{D}_{i}X^{0,\q}_{i,d_{1}}} ||v||_{\mathcal{D}_{j}X^{0,\q}_{j,d_{2}}}
\eeq

The above estimates holds true if $B(u,v)$ is replaced by $B(\bar{u},v)$ or $B(u,\bar{v})$. 

b) If $|i-j| \leq 1$, we have the following estimates:

\beq \label{be2}
||B(u,v)||_{\mathcal{D}_{k}X^{0,-\q}_{k,d_{3}}} \les 2^{\frac{k+3j}{2}} ||u||_{\mathcal{D}_{i}X^{0,\q}_{i,d_{1}}} ||v||_{\mathcal{D}_{j}X^{0,\q}_{j,d_{2}}}
\eeq

The above estimates holds true if $B(u,v)$ is replaced by $B(\bar{u},v)$ or $B(u,\bar{v})$.

\end{p3}

To simplify the exposition, we choose the work with $B(u,v)=\nabla{u} \cdot \nabla{v}$ throughout the proof. This does not restrict in any way the generality of the argument. 

\begin{proof}[Proof of Theorem \ref{pdec}]  a) We estimated $||B(u,v)||_{X_{j,d_{3}}^{0,\q}}$ for $u \in \mathcal{D}_{i}X_{i}^{0,\q,1}$ and $v \in X_{j,d_{2}}^{0,\q}$, see (\ref{bee1}). Given now the fact that $v \in \mathcal{D}_{j}X_{j,d_{2}}^{0,\q}$ we want to estimate $||B(u,v)||_{\mathcal{D}_{j}X_{j,d_{3}}^{0,\q}}$. 

One has to start with an estimate for $\chi_{Q^{m}_{j}} B(u,v)_{j,d_{3}}$ and try to commute $\chi_{Q^{m}_{j}}$ all the way next to $v$. This will be done in two steps: first commute $\chi_{Q^{m}_{j}}$ with the localization $\varphi_{j,d_{3}}$ and second with the $\nabla$. We do intend to heavily rely on the computations performed in the previous section. On the physical side, we deal with the system $Q^{m}_{j}$, while on the frequency side we deal with $A_{j,d_{3}}$ which has a $\tau$ component too and has sizes greater than the dual ones, namely $2^{-j}$, in the $\xi$ directions.  

In the same spirit with (\ref{ps1}) we can prove:

\beq \label{es1}
||\chi_{Q^{m}_{j}} h_{j,d_{3}}||_{L^{2}} \les \sum_{m'} C_{N} \langle m-m'\rangle^{-N} \sum_{|\a|+|\beta|} || \tilde{\varphi}_{j,d_{3}}^{\a} (D)  \tilde{\chi}_{Q^{m'}_{j}}^{\b} h||_{L^{2}}
\eeq

Here $\tilde{\varphi}_{j,d_{3}}^{\a}$ are a\ generalized characteristic functions of the set $A_{j,d_{3}}$ in the following sense:  $\tilde{\varphi}_{j,d_{3}}^{\a}$ is supported in $A_{j,d_{3}}$ and $||\tilde{\varphi}_{j,d_{3}}^{\a}||_{L^{\infty}} \leq C_{\a}$. 

 Matters are reduced to deal with $ \tilde{\chi}_{Q^{m}_{j}} B(u,v)$ for an arbitrary generalized characteristic function of $\tilde{\chi}_{Q^{m}_{j}}$. An exact calculus gives us:

$$ \tilde{\chi}_{Q^{m}_{j}} \nabla{v}=\nabla{(\tilde{\chi}_{Q^{m}_{j}} v)} - \nabla{(\tilde{\chi}_{Q^{m}_{j}})} v $$

We observe that it is enough to deal with the term $\nabla u \nabla{(\tilde{\chi}_{Q^{m}_{j}} v)}$. If we succeed to obtain the right estimates and then be able to sum them with respect to $m$ (over the above mentioned domain), then we will definitely be able to treat the term $\nabla{u} \cdot \nabla{(\tilde{\chi}_{Q^{m}_{j}})} v$ for the following reasons: there is no $\nabla$ on $v$ and in addition $||\nabla{\tilde{\chi}_{Q^{m}_{j}})}||_{L^{\infty}} \leq 2^{-j}$, so we are better off with a factor of $2^{-2j}$. 

The main problem we encounter in dealing with $B(u,\tilde{\chi}_{Q^{m}_{j}} v)$ is that $\tilde{\chi}_{Q^{m}_{j}}$ is not localized anymore in $A_{j,d_{2}}$ as $v$ does, which means we cannot apply directly the bilinear estimates derived before. On the other hand $ \tilde{\chi}_{Q^{m}_{j}} v$ is highly localized in $A_{j,d_{2}}$ in the following sense:

\beq \label{hard2}
||\tilde{\chi}_{Q^{m}_{j}} v||_{L^{2}_{k,d}} \les C_{P} 2^{-|k-j|P} \max{(\frac{d_{2}}{d},\frac{d}{d_{2}})}^{-P} \sum_{m'} \langle m-m'\rangle^{-P}||\chi_{Q^{m'}_{j}}v||_{L^{2}}
\eeq

We go ahead with the rest of the argument and leave the proof of this estimate for the end of the section. If we take $P \geq 3$, use the fact that $\tilde{\varphi}_{j,d_{3}}^{\a}$ is supported in $A_{j,d_{3}}$ and $||\tilde{\varphi}_{j,d_{3}}^{\a}||_{L^{\infty}} \leq C_{\a}$ and use the bilinear estimates (\ref{bee1}) we can obtain:

$$(2^{j}d_{3})^{\q}||\tilde{\varphi}_{j,d_{3}}^{\a}(D) B(u,\tilde{\chi}^{\b}_{Q^{m}_{j}} v)||_{L^{2}} \les \sum_{k,d} ||B(u, \tilde{\chi}^{\b}_{Q^{m}_{j}} v)_{k,d})||_{X_{j,d_{3}}^{0,\q}} \les$$

$$C 2^{2i} ||u||_{\mathcal{D}_{i}X^{0,\q}_{i,d_{1}}} (2^{j}d_{2})^{\q} \sum_{m'} \langle m-m'\rangle^{-P}||\chi_{Q^{m'}_{j}}v||_{L^{2}}$$

Now we can bring also the estimate in (\ref{es1}) and, if $N,P \geq 3$, we obtain:

$$\sum_{m: Q^{m}_{j} \cap L^{k}_{j} \ne \emptyset} (2^{j}d_{3})^{\q} ||\chi_{Q_{i}^{m}} B(u,v)_{j,d_{3}}||_{L^{2}} \les $$ 

$$C 2^{2i} ||u||_{\mathcal{D}_{i}X^{0,\q}_{i,d_{1}}} \sum_{m: Q^{m}_{j} \cap L^{k}_{j} \ne \emptyset} \sum_{m'} \sum_{m''}  (2^{j}d_{2})^{\q} \langle m-m'\rangle^{-N}  \langle m'-m''\rangle^{-P} ||\chi_{Q^{m''}_{j}}v||_{L^{2}} \les$$

$$C 2^{2i} ||u||_{\mathcal{D}_{i}X^{0,\q}_{i,d_{1}}} (2^{j}d_{2})^{\q} \sup_{k'} \sum_{m'':Q^{m''}_{j} \cap L^{k'}_{j} \ne \emptyset}  ||\chi_{Q^{m''}_{j}}v||_{L^{2}} \les C 2^{2i} ||u||_{\mathcal{D}_{i}X^{0,\q}_{i,d_{1}}}  ||v||_{\mathcal{D}_{j} X^{0,\q}_{j,d_{2}}} $$

Taking a $supp$ with respect to all $L^{k}_{j}$ in the aboves inequality gives us the claim in (\ref{be1}).

We owe the proof of (\ref{hard2}). For simplicity let us assume that $k=j$ and that $d, d_{2} \leq 2^{j-2}$. One can easily reproduce the argument we provide bellow for the general case. We have:

$$||\tilde{\chi}_{Q^{m}_{j}} v||^{2}_{L^{2}_{k,d}} \approx \sum_{\xi} \sum_{k} ||\varphi_{\xi,\xi^{2}+k}(D) \tilde{\chi}_{Q^{m}_{j}} v||^{2}_{L^{2}} \les$$

$$\sum_{\xi} \sum_{k} \left( \sum_{\bar{\xi}} \sum_{\bar{k}} \sum_{m'} ||\varphi_{\xi,\xi^{2}+k}(D) \tilde{\chi}_{Q^{m}_{j}} \varphi_{\bar{\xi},\bar{\xi}^{2}+\bar{k}}(D) (\chi_{Q^{m'}_{j}}v)||_{L^{2}} \right)^{2}$$

At this time we can invoke the result in (\ref{ps2}) in the following context: $Q^{m}_{j}$ is the system of cubes in the physical space of size $2^{i}$ and $A^{k}_{\xi}=\{(\eta,\tau): |\eta-\xi| \leq 2^{-i+1}, |\tau-\xi^{2}-k| \leq \q \}$ is the system of rectangles in the frequency space of size $2^{-i} \times 2^{-i} \times 1$. Since $\tilde{\chi}_{Q^{m}_{j}}$ is independent on $t$, we can ignore the $\tau$ component and then we are in the setup of the result in (\ref{ps2}), therefore:

$$||\varphi_{\xi,\xi^{2}+k}(D) \tilde{\chi}_{Q^{m}_{j}} \varphi_{\bar{\xi},\bar{\xi}^{2}+\bar{k}}(D) \chi_{Q^{m'}_{j}}v||_{L^{2}} \les $$

$$ \langle 2^{j} (\xi-\bar{\xi}) \rangle^{-P}  \langle \xi^{2}+k-\bar{\xi}^{2}-\bar{k} \rangle^{-P} \langle m-m' \rangle^{-P} \sum_{|\a|+|\b| \leq 2P} ||\tilde{\varphi}^{\a}_{\bar{\xi},\bar{\xi}^{2}+\bar{k}}(D) \tilde{\chi}^{\b}_{Q^{m'}_{j}}v||_{L^{2}}$$

The term $\langle \xi^{2}+k-\bar{\xi}^{2}-\bar{k} \rangle^{-P}$ cannot be justified via (\ref{ps2}); instead we make a simple remark: if $\langle \xi^{2}+k-\bar{\xi}^{2}-\bar{k} \rangle \geq 2$, then the actual term $\varphi_{\xi,\xi^{2}+k}(D) \tilde{\chi}_{Q^{m}_{j}} \varphi_{\bar{\xi},\bar{\xi}^{2}+\bar{k}}(D) \chi_{Q^{m'}_{j}}v$ equals $0$ since the multiplication with $\tilde{\chi_{Q_{j}^{m}}}$ does not change the $\tau$ component of the support on the Fourier side. Then we can continue with:

$$||\tilde{\chi}_{Q^{m}_{j}} v||^{2}_{L^{2}_{j,d}} \les \sum_{\xi} \sum_{k} \sum_{\bar{\xi}} \sum_{\bar{k}} \sum_{m'} \langle 2^{j} (\xi-\bar{\xi}) \rangle^{-P}  \langle \xi^{2}+k-\bar{\xi}^{2}-\bar{k} \rangle^{-P} \langle m-m' \rangle^{-P} \cdot $$

$$\sum_{|\a|+|\b| \leq 2P} ||\tilde{\varphi}^{\a}_{\bar{\xi},\bar{\xi}^{2}+\bar{k}}(D) \tilde{\chi}^{\b}_{Q^{m'}_{j}}v||^{2}_{L^{2}} \les$$

$$\max{(\frac{d_{2}}{d},\frac{d}{d_{2}})}^{-P} \sum_{|\a|+|\b| \leq 2P} \sum_{\bar{\xi}} \sum_{\bar{k}} \sum_{m'} \langle m-m' \rangle^{-P} ||\tilde{\varphi}^{\a}_{\bar{\xi},\bar{\xi}^{2}+\bar{k}}(D) \tilde{\chi}^{\b}_{Q^{m'}_{j}}v||^{2}_{L^{2}} \les$$

$$\max{(\frac{d_{2}}{d},\frac{d}{d_{2}})}^{-P} \sum_{m'} \langle m-m' \rangle^{-P} || \tilde{\chi}^{\b}_{Q^{m'}_{j}}v||^{2}_{L^{2}}$$

b) We estimated $||B(u,v)||_{X_{k,d_{3}}^{0,\q}}$ for $u \in X_{i}^{0,\q,1}$ and $v \in X_{j,d_{2}}^{0,\q}$, where $|i-j| \leq 1$, see (\ref{bee2}). Given now the fact that $v \in \mathcal{D}_{j}X_{j,d_{2}}^{0,\q}$ we want to estimate $||B(u,v)||_{\mathcal{D}_{k}X_{k,d_{3}}^{0,\q}}$. A straightforward computation gives us that:

$$||h||_{\mathcal{D}_{k}L^{2}} \les 2^{\frac{j-k}{2}} ||h||_{\mathcal{D}_{j}L^{2}}$$

This has to do with the fact that in a cube $Q^{m}_{j}$ we fit $\approx 2^{j-k}$ cubes $Q_{k}^{m'}$ on a straight line. Hence we can go ahead and estimate $||B(u,v)||_{\mathcal{D}_{j}X_{k,d_{3}}^{0,\q}}$ and bring the correction of $2^{\frac{j-k}{2}}$ at the end. Once we are in this setup we can reproduce the same argument as in part a), since $v$ comes with a $\mathcal{D}_{j}$ structure.

\end{proof}

\vspace{.1in}

\subsection{Bilinear estimates on frequency dyadic regions}
\noindent
\vspace{.1in}

In the end we want to obtain bilinear estimates on dyadic regions with respect to the frequency only.

\begin{proof}[Proof of Theorem \ref{bil}]

a) We deal first with the case when the outcome is localized at high frequency. We fix $d_{3}$ and making use of (\ref{be1}) we estimate

$$||B(u,v)||_{\mathcal{D}_{j}X^{0,-\q}_{j,d_{3}}} \les \sum_{d_{1},d_{2}} ||B(u_{\cdot,d_{1}},v_{\cdot,d_{2}})||_{\mathcal{D}_{j}X^{0,-\q}_{j,d_{3}}} \les $$

$$2^{2i} \sum_{d_{1},d_{2}} ||u_{\cdot,d_{1}}||_{\mathcal{D}_{i}X^{0,\q}_{i,d_{1}}}  ||v_{\cdot,d_{2}}||_{\mathcal{D}_{j}X^{0,\q}_{j,d_{2}}} \les 2^{2i} ||u||_{\mathcal{D}_{i}X^{0,\q,1}_{i}}  ||v||_{\mathcal{D}_{j}X^{0,\q,1}_{j}}$$

Summing up with $d_{3}$ and passing to general $s$ gives us (\ref{be11}).

In the case when $|i-j| \leq 1$ and $k \leq j-5$, we estimate in the same way, this time making use of (\ref{be2}), to obtain (\ref{be22}).  

b) We decompose 

$$v_{\cdot, \geq 2^{-i}}=\sum_{2^{-i} \leq d' \leq 2^{i}} v_{\cdot, d'} + \sum_{d_{2} \geq 2^{i+1}} v_{\cdot, d_{2}}$$ 

and notice that $ \hat{u} * \sum_{2^{-i} \leq d' \leq 2^{i}} \hat{v}_{\cdot, d'}$ is essentially localized at distance less than $2^{i}$ from $P$ while $\hat{u} * \hat{v}_{\cdot,d_{2}}$ is localized essentially at distance $d_{2}$ from $P$ for any $d_{2} \geq 2^{i+1}$. This happens because $\hat{u}$ is localized at frequency $2^{i}$.

We fix $d_{3} \geq 2^{-i}$ and as in part a) we estimate:

$$||B(u, \sum_{2^{-i} \leq d' \leq 2^{i}} v_{\cdot, d'})||_{X^{0,-\q}_{j,d_{3}}} \les 2^{2i}  ||u||_{X^{0,\q,1}_{i}} ||\sum_{2^{-i} \leq d' \leq 2^{i}} v_{\cdot, d'}||_{X^{0,\q,1}_{j}}$$

In a similar manner we can conclude that for any $d_{2} \geq 2^{i}$ we obtain:

$$||B(u, v_{\cdot, d_{2}})||_{X^{0,-\q}_{j,d_{2}}}  \les 2^{2i}  ||u||_{X^{0,\q,1}_{i}} || v_{\cdot, d_{2}}||_{X^{0,\q,1}_{j}}$$

Taking into account the above observation above the localization of the interactions, we sum up with respect to $d_{3}$, for $2^{-i} \leq d_{3} \leq 2^{i}$ and then with respect to $d_{2} \geq 2^{i}$, to obtain:

$$||B(u, v)||_{X^{0,-\q,1}_{j, \leq 2^{-i}}}  \les 2^{2i} i ||u||_{X^{0,\q,1}_{i}} ||v||_{X^{0,\q,1}_{j, \leq 2^{-i}}} $$

Passing to general $s$ gives us (\ref{be22}).

\end{proof}

\vspace{.1in}

\section{Bilinear estimates involving the {\mathversion{bold}$Y$} spaces}

\vspace{.1in}

In the previous section we have just seen that the theory of bilinear estimates cannot be completely closed in the $X^{s,\q,1}$ spaces. This is the reason for introducing a more refined structure to measure our solutions, namely the wave-packet one. We concluded that the interactions causing problems in the $X^{s,\q,1}$ theory are the low-high ones. This is why we need to complete Theorem \ref{bil} with a result for this particular case.

\begin{t1} \label{tb2}

Assume we have $5i \leq j$. We have the bilinear estimates:

\beq \label{b7}
||B(u,v)||_{\mathcal{D}W^{s}_{j}} \les i 2^{(2-s)i} ||u||_{\mathcal{D}Z^{s}_{i}} ||v||_{\mathcal{D}Z^{s}_{j}}
\eeq

The estimate remains valid if $B(u,v)$ is replaced by $B(\bar{u},v)$ or $B(u,\bar{v})$. 

\end{t1}

 In what follows we make few important remarks for the rest of this section. The first one comes from the hypothesis of our theorem. 

\begin{r3} 
We work under the hypothesis that $5i \leq j$. 
\end{r3}
The result in (\ref{be22}) shows that it is fine to use the $X^{s,\q,1}$ structure to measure the low frequency and part of the high frequency (both input and output) at distance greater than $2^{-i}$ from $P$. Thus we shall obtain estimates for:

\beq \label{t1}
X_{i}^{0,\q,1} \cdot Y_{j, \leq 2^{-i}} \rightarrow \mathcal{Y}_{j, \leq 2^{-i}} \ + \ X^{0,\q,1}_{j,\geq 2^{-i}}; \  X_{i}^{0,\q,1} \cdot X^{0,\q,1}_{j, \geq 2^{-i}} \rightarrow \mathcal{Y}_{j, \leq 2^{-i}} 
\eeq

 We also need the corresponding estimates when we involve conjugates of these spaces. The condition $5i \leq j$ implies that the the low frequency does not see the curvature of the parabola at the high frequency, in other words the parabola at high frequency is flat in these interactions. This is why the estimates for $B(\bar{u}_{i},v_{j})$ are similar to the ones for $B(u_{i},v_{j})$. 
 
 If we have to deal with $B(u_{i},\bar{v}_{j})$, a simple geometric argument shows that the interaction is localized at high frequency and in a region with $\tau \leq 0$. This makes these estimates weaker than the ones in (\ref{t1}).

\begin{r3}
Once we get one of the estimates in (\ref{t1}), we trivially get the corresponding ones with conjugate spaces. 
\end{r3}

We have to involve and recover decay in these estimates. We prove:

$$||B(u,v)||_{\mathcal{R}W^{s}_{j}} \les i^{\frac{3}{2}} 2^{(1-s)i} ||u||_{\mathcal{R}\mathcal{D}Z^{s}_{i}} ||v||_{\mathcal{R}Z^{s}_{j}}$$

\noindent
and the similar ones. In the end we obtain the estimates with decay on all terms by a similar argument as in section \ref{cd}.  

\begin{r3}

We first prove the estimates without involving decay on the bilinear term and on the high frequency. But we do involve decay on the low frequency. 

\end{r3}

These being said, we can start the preparations for this section.

\vspace{.1in}

\subsection{Basic estimates}

\noindent
\vspace{.1in}

This section is concerned with providing results of type $Y \cdot \mathcal{D}L^{2} \rightarrow \mathcal{Y}$, $Y \cdot \mathcal{D} L^{2} \rightarrow L^{2}$ and $L^{2} \cdot \mathcal{D}L^{2} \rightarrow \mathcal{Y}$.

\begin{l5}

Let $g \in L^{2}$ such that $\hat{g}$ is supported in a tube of size $2^{-i} \times 2^{-i} \times 1$. We have the estimate:

\beq \label{c0}
\sum_{m} ||g||^{2}_{L^{\infty}(Q^{m,l}_{i})} \les 2^{-2i} ||g||^{2}_{L^{2}}
\eeq

\end{l5}

\begin{proof}

The support of $\hat{g}$ is a tube with volume $2^{-2i}$ therefore we have:

$$||g||_{L^{\infty}(Q^{m,l}_{i}))} \les 2^{-i} \sum_{(m',l') \in Z^{3}} C_{N} \langle (m,l)-(m',l') \rangle^{-N} ||g||_{L^{2}(Q^{m',l'}_{i})}$$

If we chose $N \geq 4$, then we use Cauchy-Schwartz and get:

$$||g||^{2}_{L^{\infty}(Q^{m,l}_{i}))} \les  2^{-2i} \sum_{(m',l') \in \Z^{3}} C^{2}_{N} \langle (m,l)-(m',l') \rangle^{-N} ||g||^{2}_{L^{2}(Q^{m',l'}_{i})} $$

We can perform the summation with respect to $(m,l)$:

$$\sum_{(m,l)} ||g||^{2}_{L^{\infty}(Q^{m,l}_{i}))} \les  2^{-2i} \sum_{m,l} \sum_{m',l'}   \langle (m,l)-(m',l') \rangle^{-N} ||g||^{2}_{L^{2}(Q^{m',l'}_{i} )} \les 2^{-2i} ||g||^{2}_{L^{2}}$$

In the last line we use again the fact that if $N \geq 4$, then we have:

$$\sum_{m,l}  \langle (m,l)-(m',l') \rangle^{-N} \les 1$$

This is enough to justify the claim.

\end{proof}

\begin{l5}

Let $\hat{g}$ be supported in $A_{i,d}$ where $d \leq 2^{i-2}$. For any $p \in \Z$ we have:

\beq \label{dd1}
\sum_{\xi} \sum_{k=2^{i-1}d}^{2^{i+1}d} ||\chi_{Q_{i}^{m}} g_{\xi,\xi^{2}+k}||^{2}_{L^{\infty}} \les C_{N} 2^{-i} \sum_{m'} \langle m-m' \rangle^{-N} ||\chi_{Q^{m'}_{i}} g||^{2}_{L^{2}}
\eeq

\beq \label{dd2}
\sum_{m: Q^{m}_{i} \cap L^{p}_i \ne \emptyset} \left( \sum_{\xi} \sum_{k=2^{i-1}d}^{2^{i+1}d} ||\chi_{Q_{i}^{m}} g_{\xi,\xi^{2}+k}||^{2}_{L^{\infty}}  \right)^{\q} \les 2^{-i} ||g||_{\mathcal{D}_{i}L^2}
\eeq

\end{l5}

\begin{proof}

The support of $\hat{g}_{\xi,\xi^{2}+k}$ is $2^{-i} \times 2^{-i} \times 1$, hence:

$$||\chi_{Q_{i}^{m}}g_{\xi,\xi^{2}+k}||_{L^{\infty}} \les C_{N} 2^{-i} \sum_{m'} \langle m-m' \rangle^{-N} ||\chi_{Q^{m'}_{i}} g_{\xi,\xi^{2}+k}||_{L^{2}}$$

Then, (\ref{dd1}) amounts to proving:

\beq \label{dd3}
\sum_{\xi} \sum_{k} ||\chi_{Q_{i}^{m}} g_{\xi,\xi^{2}+k}||^{2}_{L^{2}} \les C_{N} \sum_{m'} \langle m-m' \rangle^{-N} ||\chi_{Q^{m'}_{i}} g||_{L^{2}}
\eeq

For fixed $\xi$ we have the obvious:

$$\sum_{k} ||\chi_{Q_{i}^{m}} g_{\xi,\xi^{2}+k}||^{2}_{L^{2}} \approx ||\chi_{Q_{i}^{m}} \sum_{k} g_{\xi,\xi^{2}+k}||^{2}_{L^{2}}$$ 

\noindent
since $\chi_{Q_{i}^{m}}$ is a cut in the $x$ space while $(\varphi_{\xi,\xi^{2}+k})_{k}$ is a cut in the $\tau$ direction. Hence it is enough to prove (\ref{dd3}) in the particular case $d=2^{-i}$(i.e. $k=0$):

$$\sum_{\xi} ||\chi_{Q_{i}^{m}} g_{\xi,\xi^{2}}||^{2}_{L^{2}} \les C_{N} \sum_{m'} \langle m-m' \rangle^{-N} ||\chi_{Q^{m'}_{i}} g||_{L^{2}}$$

We can write:

$$\chi_{Q^{m}_{i}}g_{\xi,\xi^{2}}=\chi_{Q^{m}_{i}} (\sum_{m'}  \chi_{Q^{m'}_{i}} g)_{\xi,\xi^{2}}$$

Invoking the results (and notations) from section \ref{ar}, see also the adjustments in section \ref{cd}, we claim:

\beq \label{dd4}
||\chi_{Q^{m}_{i}} ( \chi_{Q^{m'}_{i}} g)_{\xi,\xi^{2}}||_{L^{2}} \les C_{N} \langle m-m'\rangle^{-N} \sum_{|\a|+|\b| \leq N} || \tilde{\varphi}^{\a}_{\xi,\xi^{2}}(D) \tilde{\chi}^{\b}_{Q^{m'}_{i}} g)||_{L^{2}}
\eeq

Then we sum with respect to $m'$ and use Cauchy-Schwartz to obtain:

$$||\chi_{Q^{m}_{i}} g_{\xi,\xi^{2}}||^{2}_{L^{2}} \les C_{N} \langle m-m'\rangle^{-N} \sum_{|\a|+|\b| \leq N} || \tilde{\varphi}^{\a}_{\xi,\xi^{2}}(D) \tilde{\chi}^{\b}_{Q^{m'}_{i}} g)||^{2}_{L^{2}}$$ 

Recalling (\ref{ps5}), both in space and frequency, we continue with:

$$\sum_{\xi} ||\chi_{Q^{m}_{i}} g_{\xi,\xi^{2}}||^{2}_{L^{2}} \les C_{N} \langle m-m'\rangle^{-N} \sum_{|\a|+|\b| \leq N} \sum_{\xi} || \tilde{\varphi}^{\a}_{\xi,\xi^{2}}(D) \tilde{\chi}^{\b}_{Q^{m'}_{i}} g)||^{2}_{L^{2}} \les$$

$$C_{N} \langle m-m'\rangle^{-N} \sum_{\a} ||\chi_{Q^{m'}_{i}} g||^{2}_{L^{2}}$$

Passing from (\ref{dd1}) to (\ref{dd2}) is a matter of algebraic computations.

\end{proof}

In a similar way we can prove the following result:

\begin{l5}

If $\hat{g}$ is supported in $A_{i,d}$ for $d \geq 2^{i-2}$, then for any $p \in \Z$:

\beq \label{dd6}
\sum_{m: Q^{m}_{i} \cap L^{p}_i \ne \emptyset} \left( \sum_{\xi} \sum_{l=2^{2i-2}d}^{2^{2i+2}} ||\chi_{Q_{i}^{m}} g_{\xi,l}||^{2}_{L^{\infty}}  \right)^{\q} \les 2^{-i} ||g||_{\mathcal{D}_{i}L^2}
\eeq

\end{l5}

\vspace{.1in}

For each $\a \in Z^{2}$ we define $A_{\a}=\{m \in \Z^{2}: m_{1} \in [2^{j}(2\a-1),2^{j}(2\a+1)], m_{2} \in [2^{i}(2\a-1),2^{i}(2\a+1)]\}$. We have the following result:

\begin{l5} \label{srty}

The families $(T^{m,l}_{\eta})_{m \in A_{\a}}$ and $(T^{m,l}_{\eta + \xi})_{m \in A_{\b}}$ contain disjoint tubes unless $|\a-\b|=\max(|\a_{1}-\b_{1}|,|\a_{2}-\b_{2}|) \leq 2$; in other words if $T^{m,l}_{\eta + \xi} \cap T^{m',l}_{\eta + \xi} \ne \emptyset$, where  $m \in A_{\a}$ and $m' \in A_{\b}$, then $|\a-\b| \leq 2$.

\end{l5}

\begin{proof}

It is enough to prove the result in the case $l=0$. Let us assume that there is $(x,t) \in T^{m,0}_{\eta + \xi} \cap T^{m',0}_{\eta + \xi}$, where  $m \in A_{\a}$ and $m' \in A_{\b}$. Then:

$$||x-m+t \eta|| \leq \sqrt{2} \ \ \ \mbox{and} \ \ \ |x-m'+t(\xi+\eta)| \leq \sqrt{2}$$

\noindent
which implies $||m-m'+t\xi|| \leq 2\sqrt{2}$. Recalling that $t \in [0,1]$, $||\xi|| \approx 2^{i}$, $i < j$ and the definition of $A_{\a}, A_{\beta}$ we obtain the claim.

\end{proof}

\begin{l5} \label{lavt}

For each $m,m'$ there is essentially only one $m''$ such that $Q_{i}^{m,l} \cap T_{\eta}^{m',l} \cap T_{\eta+\xi}^{m'',l} \ne \emptyset$; more precisely, there are at most $5$ $m''$'s with this property.

\end{l5}

\begin{proof}

The underlying idea is that the intersection $T_{\eta}^{m',l} \cap T_{\eta+\xi}^{m'',l}$ is a subtube of sizes $2^{j-i}$ in the long direction and $2^{j-i} \geq 2^{i}$, the later being the size of the cube $Q^{m,l}_{i}$. One can formalize an explicit proof.

\end{proof}

For each $\a$, we define by $B_{\a}=\{m: Q_{i}^{m,l} \cap T_{\eta}^{m',l} \ne \emptyset \ \ \mbox{for} \ m' \in A_{\a} \}$. Notice that the family of tubes $(T_{\eta}^{m',l})_{m' \in A_{\a}}$ fill up a parallelepiped of sizes $2^{j+1} \times 2^{i+1} \times 1$ (last one in the $t$ direction) and the longest side is in the direction of $\eta$. Hence if $L$ is the line in $\R^{2}$ passing through the origin in the direction of $\eta$, then there is a $k \in \Z$ such that: 

\beq \label{lkjh}
B_{\a} \subset \{m: Q_{i}^{m} \cap (L^{k-1}_{i} \cup L_{i}^{k} \cup L^{k+1}_{i}) \ne \emptyset\}
\eeq

We conclude with the main result of this section.

\begin{l5} We have the estimate:

\beq \label{m9}
|| f \cdot  g||_{L^{2}} \les  2^{-\frac{i+j}{2}} ||f||_{Y_{j}} ||g||_{L^{2}}
\eeq

\end{l5}

\begin{proof}

 For $m' \in A_{\a}$, we have:

$$||f \cdot g||^{2}_{L^{2}(T^{m',l}_{\eta})} \les \sum_{m \in B_{\a}} ||f \cdot g||^{2}_{L^{2}(Q_{i}^{m,l} \cap T^{m',l}_{\eta})} \les$$

$$ 2^{i-j} \sum_{m \in B_{\a}} ||f||^{2}_{L_{t}^{\infty}L^{2}_{x}(T^{m',l}_{\eta})} ||g||^{2}_{L^{\infty}(Q_{i}^{m,l})} \les $$

$$2^{i-j} ||f||^{2}_{L_{t}^{\infty}L^{2}_{x}(T^{m',l}_{\eta})} \sum_{m \in B_{\a}} ||g||^{2}_{L^{\infty}(Q_{i}^{m,l})} \les   2^{-i-j} ||f||^{2}_{L_{t}^{\infty}L^{2}_{x}(T^{m',l}_{\eta})} ||g||^{2}_{L^{2}}$$

In the last line we have used the result in (\ref{c0}). We sum the above estimate with respect to $(m,l)$ over $Z^{3}$ to obtain (\ref{m9}).

\end{proof}

The next Lemma is a geometrical one. We work with $f=f_{\eta, \leq 2^{-i}}$ and $g=g_{\xi^{0},l}$, $\xi^{0} \in \Xi^{i}, l \in \Z$ where $|\eta| \approx 2^{j}$ and $|(\xi^{0},l)| \approx 2^{i}$.

\begin{l5}
Assume $d \geq 2^{-i}$. If $\hat{g} * \hat{f}$ is supported in a region where $|\tau-\xi^{2}| \leq d$ then $|\cos{\a}| \leq |\xi^{0}|^{-1}d $, 
where $\a$ is the angle between $\xi^{0}$ and $\eta$. 

\end{l5}

\begin{proof}

 $\hat{f}$ is supported in a region where $|\tau_{2}-\eta^{2}| \leq 2^{-i}|\eta|$, while $\hat{g}$ is supported in a region where $|\xi-\xi^{0}| \leq 2^{-i}$ and $|\tau_{1}-l| \leq \q$. A generic point in the support of $\hat{f}*\hat{g}$ is of type $(\xi_{1}+\xi_{2},\tau_{1}+\tau_{2})$ where $(\xi_{1},\tau_{1})$ is in the support of $\hat{g}$ and $(\xi_{2},\tau_{2})$ is in the support of $\hat{f}$. We want this point to satisfy $|\tau_{1}+\tau_{2}-(\xi_{1}+\xi_{2})^{2}| \leq d$.

We have $|\tau_{1} -\xi_{1}^{2}| \leq 2^{2i} \leq 2^{j-i}$, $\D |\xi_{1}| \approx 2^{-i}$, $\D |\eta| \approx 1$, therefore the condition is equivalent to $|2\eta \cdot \xi^{0}| \leq 2^{j}d$. This implies the conclusion of the Lemma.

\end{proof}

\begin{l5} \label{g1}
 For fixed $\xi$ and $k$, the interactions $\hat{g}_{\xi,\xi^{2}+k}*\hat{f}_{\eta,\leq 2^{-i}}$ have disjoint supports with respect to $\eta$; same is true for $\hat{g}_{\xi,l}*f_{\eta,\leq 2^{-i}}$.

\end{l5}

\begin{proof}
The sizes of the support of $\hat{g}_{\xi,\xi^{2}+k}$ are $2^{-i} \times 2^{-i} \times 1$. The support of $\hat{v}_{\eta, \leq 2^{-i}}$ is a parallelepiped of sizes $2^{-i} \times 1 \times 2^{j}$ whose longest side is tangent to $P$. The key property is that we can translate the support of $\hat{u}_{\xi, \leq 2^{-i}}$ so that it is included in the support of $\hat{v}_{\eta, \leq 2^{-i}}$ (by simply translating the center of the first to the center of the second). Therefore the support of $\hat{v}_{\eta, \leq 2^{-i}} * \hat{u}_{\xi, \leq 2^{-i}}$ is a translate of the support of $\hat{v}_{\eta, \leq 2^{-i}}$ by the vector $(\xi, \xi^{2})$. Therefore if we keep $\xi$ and $k$ fixed and take  $\eta \ne \eta'$  both in $A_{\xi}$, then the supports of $\hat{v}_{\eta, \leq 2^{-i}} * \hat{u}_{\xi, \leq 2^{-i}}$ and $\hat{v}_{\eta', \leq 2^{-i}} * \hat{u}_{\xi, \leq 2^{-i}}$ are disjoint.

\end{proof}

\vspace{.1in}

\subsection{Estimates: {\mathversion{bold}$ \mathcal{D}X^{0,\q,1}_{i} \cdot Y_{j, \leq 2^{-i}} \rightarrow \mathcal{Y}_{j, \leq 2^{-i}}$}}

\noindent
\vspace{.1in}

The main result of this section is the following:

\begin{p3} We have the estimate:

\beq \label{a0}
|| v_{j, \leq 2^{-i}} \cdot  u_{i}||_{\mathcal{Y}_{j, \leq 2^{-i}}} \les 2^{i-j} ||v_{j, \leq 2^{-i}}||_{Y_{j}} \cdot ||u_{i}||_{\mathcal{D}X^{0,\q,1}}
\eeq

\end{p3}

This result is a direct consequence of the following estimates:

\beq \label{a1}
|| f_{j, \leq 2^{-i}} \cdot  g_{i, \leq 2^{i-2}}||_{\mathcal{Y}_{j, \leq 2^{-i}}} \les 2^{i-j} ||f_{j, \leq 2^{-i}}||_{Y_{j}} \cdot ||g_{i, \leq 2^{i-2}}||_{\mathcal{D}X^{0,\q,1}}
\eeq

\beq \label{a2}
|| f_{j, \leq 2^{-i}} \cdot  g_{i, \geq 2^{i-2}}||_{\mathcal{Y}_{j, \leq 2^{-i}}} \les 2^{i-j} ||f_{j, \leq 2^{-i}}||_{Y_{j}} \cdot ||g_{i, \geq 2^{i-2}}||_{\mathcal{D}X^{0,\q,1}}
\eeq

\beq \label{a3}
|| f_{j, \leq 2^{-i}} \cdot  g_{i}||_{\mathcal{Y}_{j, \leq 2^{-i}}} \les 2^{i-j} ||f_{j, \leq 2^{-i}}||_{Y_{j}} \cdot ||g_{i}||_{\mathcal{D}X^{0,\q,1}}
\eeq

\beq \label{a4}
|| f_{j, \leq 2^{-i}} \cdot  g_{i}||_{\mathcal{D}\mathcal{Y}_{j, \leq 2^{-i}}} \les 2^{i-j} ||f_{j, \leq 2^{-i}}||_{\mathcal{D}Y_{j}} \cdot ||g_{i}||_{\mathcal{D}X^{0,\q,1}}
\eeq

\begin{proof} Throughout this section we use the following decompositions:

\beq \label{dec1}
f_{j, \leq 2^{-i}}= \sum_{\eta \in \Xi} f_{\eta, \leq 2^{-i}}
\eeq

\beq \label{dec2}
g_{i}=g_{i, \leq 2^{i-2}} + g_{i, \geq 2^{i-2}} = \sum_{k} \sum_{\xi \in \Xi^{i}} g_{\xi, \xi^{2} \pm k}+  \sum_{\xi \in \Xi^{i}} \sum_{l \in I_{\xi}} g_{\xi,l}
\eeq

For more details about the decomposition in (\ref{dec2}), see (\ref{aux2}) and (\ref{aux1}). We do not want to bother about carrying the $\pm$ in $g_{\xi, \xi^{2} \pm k}$ in all computations. We choose to work only with the $g_{\xi, \xi^{2} + k}$ ($k$ will be positive) and completing the argument for both choices of sign is a trivial matter. 

We prove first (\ref{a1}). We make use of the decompositions in (\ref{dec1}) and (\ref{dec2}). We define $\Theta_{p}=\{\eta \in \Xi, |\eta| \approx 2^{j}: \arg{\eta} \in [(p-\q)2^{i-j}\frac{\pi}{2},(p+\q)2^{i-j}\frac{\pi}{2}]\}$. The size of $\sum_{\eta \in \Theta_{p}} \hat{f}_{\eta, \leq 2^{-i}}$ in the angular direction is $\approx 2^{i}$, therefore the interactions $\hat{g}_{i} * \sum_{\eta \in \Theta_{p}} \hat{f}_{\eta, \leq 2^{-i}}$ have disjoint support with respect to $p$. As a consequence:

\beq \label{k2}
||S_{j, \leq 2^{-i}}( g_{i} \cdot f_{j, \leq 2^{-i}})||^{2}_{\mathcal{Y}} \approx \sum_{p} ||S_{j, \leq 2^{-i}}(g_{i} \cdot \sum_{\eta \in \Theta_{p}} f_{\eta, \leq 2^{-i}})||^{2}_{\mathcal{Y}}
\eeq

We decompose:

$$S_{j, \leq 2^{-i}} (g_{i,\leq 2^{i-2}} \cdot \sum_{\eta \in \Theta_{p}} f_{\eta, \leq 2^{-i}})= S_{j, \leq 2^{-i}} \left( \sum_{d \leq 2^{i-2}} \sum_{\xi} \sum_{k=2^{i-1}d}^{2^{i+1}d} g_{\xi, \xi^{2}+k} \right) \left( \sum_{\eta \in \Theta_{p}} f_{\eta, \leq 2^{-i}} \right) $$

From Lemma \ref{g1} we know that $g_{\xi, \xi^{2}+k} \cdot f_{\eta, \leq 2^{-i}}$ is supported in $A_{j,\leq 2^{-i}}$ iff $|\cos{\a}| \leq 2^{-2i}$, where $\a$ is the angle between $\xi$ and $\eta$. The angle between any two $\eta$'s in $\Theta_{p}$ is at most $2^{i-j} \leq 2^{-2i}$ and the angle between any two $\xi$'s in $\Xi^{i}$ is either at least $2^{-2i}$ or the same. Therefore all $\xi's$ involved in the above summation have the same angular localization; we just keep this in mind and not formalize it. What is important is that we sum over a set containing $\approx 2^{2i}$ $\xi$'s. We continue with:  

$$S_{j, \leq 2^{-i}} \left( \sum_{\xi} \sum_{k} g_{\xi, \xi^{2}+k} \right) \left( \sum_{\eta \in \Theta_{p}} f_{\eta, \leq 2^{-i}} \right)=$$

$$ \sum_{k} \sum_{\xi} \sum_{\eta \in \Theta_{p}} (g_{\xi, \xi^{2}+k} \cdot f_{\eta, \leq 2^{-i}})_{\xi+\eta, \leq 2^{-i}}=$$

$$\sum_{\a} \sum_{l} \sum_{m \in B^{\a}} \sum_{m' \in A^{\a}} \sum_{k} \sum_{\xi} \sum_{\eta \in \Theta_{p}} (\chi_{Q_{i}^{m}} g_{\xi, \xi^{2}+k} \cdot \chi_{T_{\eta}^{m',l}} f_{\eta, \leq 2^{-i}})_{\xi+\eta, \leq 2^{-i}}$$  

For fixed $\a$ and $l$, $\sum_{m \in B^{\a}} \chi_{T_{\eta}^{m',l}} f_{\eta, \leq 2^{-i}}$ is essentially supported (in the physical space) in a parallelepiped of sizes $2^{i} \times 2^{j} \times 1$ which is independent of $\eta \in \Theta_{p}$. The position of this parallelepiped is function of $\a$ and $l$. Hence we have:

\beq \label{k1}
||S_{j, \leq 2^{-i}} \left( \sum_{\xi} \sum_{k} g_{\xi, \xi^{2}+k} \right) \left( \sum_{\eta \in \Theta_{p}} f_{\eta, \leq 2^{-i}} \right)||^{2}_{\mathcal{Y}} \approx
\eeq

$$\sum_{\a} \sum_{l} ||\sum_{m \in B^{\a}} \sum_{m' \in A^{\a}} \sum_{k} \sum_{\xi} \sum_{\eta \in \Theta_{p}} (\chi_{Q_{i}^{m}} g_{\xi, \xi^{2}+k} \cdot \chi_{T_{\eta}^{m',l}} f_{\eta, \leq 2^{-i}})_{\xi+\eta, \leq 2^{-i}}||^{2}_{\mathcal{Y}}$$  

We fix $\a$, $l$ and $m \in B_{\a}$. Without losing the generality of the argument, we choose $l=0$. We fix $\xi $ and $k \in [2^{i-1}d,2^{i+1}d]$. We also want to drop the notation relocalization $(\cdot)_{\eta+\xi,\leq 2^{-i}}$ and we can do that by making the convention that $\chi g_{\xi,\xi^{2}+k} \cdot \chi f_{\eta}$ has to be measured in $\mathcal{Y}_{\eta+\xi}$.

We continue with:

$$||\sum_{m' \in A^{\a}} \sum_{\eta} \chi_{Q_{i}^{m}} g_{\xi, \xi^{2}+k} \cdot \chi_{T_{\eta}^{m',0}} f_{\eta, \leq 2^{-i}}||_{\mathcal{Y}}^{2} \approx$$

$$\sum_{m' \in A^{\a}} \sum_{\eta \in \Theta_{p}} \sum_{m''} ||\chi_{Q_{i}^{m}} g_{\xi, \xi^{2}+k} \cdot \chi_{T_{\eta}^{m',0}} f_{\eta, \leq 2^{-i}}||^{2}_{L^{1}_{t}L^{2}_{x}(T^{m'',0}_{\eta+\xi})}$$

For fixed $m'$, let $m''$ be such that $Q_{i}^{m} \cap T_{\eta}^{m',0} \cap T_{\eta+\xi}^{m'',0} \ne \emptyset$. The size of this intersection in the direction of $t$ is $\approx 2^{i-j}$, therefore we can estimate:

$$||\chi_{Q_{i}^{m}} g_{\xi, \xi^{2}+k} \cdot \chi_{T^{m',0}_{\eta}}f||_{L^{1}_{t}L^{2}_{x}(T_{\eta+\xi}^{m'',0})} \les 2^{i-j} ||\chi_{Q_{i}^{m}} g_{\xi, \xi^{2}+k} \cdot \chi_{T^{m',0}_{\eta}}f_{\eta, \leq 2^{-i}}||_{L^{\infty}_{t}L^{2}_{x}} \les$$

$$2^{i-j} ||\chi_{Q_{i}^{m}} g_{\xi, \xi^{2}+k}||_{L^{\infty}} || \chi_{T^{m',0}_{\eta}}f_{\eta, \leq 2^{-i}}||_{L^{\infty}_{t}L^{2}_{x}}$$

Taking into account the result of Lemma \ref{lavt} we can perform the $l^{2}_{m'}$ summation and obtain:

$$||\chi_{Q_{i}^{m}} g_{\xi, \xi^{2}+k} \cdot \sum_{m' \in A_{\a}} \chi_{T^{m',0}_{\eta}}f_{\eta, \leq 2^{-i}}||_{\mathcal{Y}_{\xi+\eta}} \les $$

$$2^{i-j} ||\chi_{Q_{i}^{m}} g_{\xi, \xi^{2}+k}||_{L^{\infty}} ||\sum_{m' \in A_{\a}} \chi_{T^{m',0}_{\eta}}f_{\eta, \leq 2^{-i}}||_{Y_{\eta}}$$

Next we can perform the $l^{2}_{\eta}$ summation to obtain:

$$||\chi_{Q_{i}^{m}} g_{\xi, \xi^{2}+k} \cdot \sum_{\eta \in \Theta_{p}} \sum_{m' \in A_{\a}} \chi_{T^{m',0}_{\eta}}f_{\eta, \leq 2^{-i}}||_{\mathcal{Y}_{\xi+\eta}} \les$$

$$ 2^{i-j} ||\chi_{Q_{i}^{m}} g_{\xi, \xi^{2}+k}||_{L^{\infty}} \left( \sum_{\eta \in \Theta_{p}} ||\sum_{m' \in A_{\a}} \chi_{T^{m',0}_{\eta}}f_{\eta, \leq 2^{-i}}||^{2}_{Y_{\eta}} \right)^{\q}$$

We fix $d \in I_{i}$ and perform the summation with respect to $\xi$ and $k \in [2^{i-1}d, 2^{i+1}d]$:

$$||\sum_{\xi} \sum_{k} \chi_{Q_{i}^{m}} g_{\xi, \xi^{2}+k} \cdot \sum_{\eta \in \Theta_{p}} \sum_{m' \in A_{\a}} \chi_{T^{m',0}_{\eta}}f_{\eta, \leq 2^{-i}}||_{\mathcal{Y}_{\xi+\eta}} \les $$

$$2^{i-j} 2^{i}(2^{i}d)^{\q} \left( \sum_{\xi} \sum_{k}||\chi_{Q_{i}^{m}} g_{\xi, \xi^{2}+k}||^{2}_{L^{\infty}} \right)^{\q} \left( \sum_{\eta \in \Theta_{p}} ||\sum_{m' \in A_{\a}} \chi_{T^{m',0}_{\eta}}f_{\eta, \leq 2^{-i}}||^{2}_{Y_{\eta}} \right)^{\q}$$

We sum up with respect to $m \in B_{\a}$:

$$||\sum_{m \in B_{\a}}\sum_{\xi} \sum_{k} \chi_{Q_{i}^{m}} g_{\xi, \xi^{2}+k} \cdot \sum_{\eta \in \Theta_{p}} \sum_{m' \in A_{\a}} \chi_{T^{m',0}_{\eta}}f_{\eta, \leq 2^{-i}}||_{\mathcal{Y}_{\xi+\eta}} \les $$

$$2^{i-j} 2^{i}(2^{i}d)^{\q}\sum_{m \in B_{\a}} \left( \sum_{\xi} \sum_{k}|| \chi_{Q_{i}^{m}} g_{\xi, \xi^{2}+k}||^{2}_{L^{\infty}} \right)^{\q} \left( \sum_{\eta \in \Theta_{p}} ||\sum_{m' \in A_{\a}} \chi_{T^{m',0}_{\eta}}f_{\eta, \leq 2^{-i}}||^{2}_{Y_{\eta}} \right)^{\q} \les$$

$$2^{i-j}  ||g_{i,d}||_{\mathcal{D}X^{0,\q}} \left( \sum_{\eta \in \Theta_{p}} ||\sum_{m' \in A_{\a}} \chi_{T^{m',0}_{\eta}}f_{\eta, \leq 2^{-i}}||^{2}_{Y_{\eta}} \right)^{\q} $$

In the last inequality we have used (\ref{dd2}). We sum with respect to $d$:

$$||\sum_{m \in B_{\a}}\sum_{\xi} \sum_{k} \chi_{Q_{i}^{m}} g_{\xi, \xi^{2}+k} \cdot \sum_{\eta \in \Theta_{p}} \sum_{m' \in A_{\a}} \chi_{T^{m',0}_{\eta}}f_{\eta, \leq 2^{-i}}||_{\mathcal{Y}_{\xi+\eta}} \les $$

$$2^{i-j} ||g_{i,\leq 2^{i-2}}||_{\mathcal{D}X^{0,\q,1}} \left( \sum_{\eta \in \Theta_{p}} ||\sum_{m' \in A_{\a}} \chi_{T^{m',0}_{\eta}}f_{\eta, \leq 2^{-i}}||^{2}_{Y_{\eta}} \right)^{\q} $$

Now we make use of (\ref{k1}) and sum up with respect to $\a$ and $l$ to obtain:

$$||S_{j, \leq 2^{-i}} \left( \sum_{\xi} \sum_{k} g_{\xi, \xi^{2}+k} \right) \left( \sum_{\eta \in \Theta_{p}} f_{\eta, \leq 2^{-i}} \right)||^{2}_{\mathcal{Y}} \les $$

$$2^{i-j}  ||g_{i,\leq 2^{i-2}}||_{\mathcal{D}X^{0,\q,1}} \left( \sum_{\eta \in \Theta_{p}} ||f_{\eta, \leq 2^{-i}}||^{2}_{Y_{\eta}} \right)^{\q} $$

In the end we use (\ref{k2}) to perform the summation with respect to $p$ and obtain the claim in (\ref{a1}).

The argument for (\ref{a2}) is carried on in the same fashion. We have the estimate (\ref{dd6}) to replace (\ref{dd2}) in this case.

(\ref{a3}) is the sum of (\ref{a1}) and (\ref{a2}).

(\ref{a4}) is the sum of the variants of  (\ref{a1}) and (\ref{a2}) with decay. We sketch the proof for: 

\beq \label{a5}
|| v_{j, \leq 2^{-i}} \cdot  u_{i, \leq 2^{i-2}}||_{\mathcal{D}\mathcal{Y}_{j, \leq 2^{-i}}} \les 2^{i-j} ||v_{j, \leq 2^{-i}}||_{\mathcal{D}Y_{j}} \cdot ||u_{i, \leq 2^{i-2}}||_{\mathcal{D}X^{0,\q,1}}
\eeq

We follow the steps in the proof of (\ref{a1}). For fixed $\tilde{m}$ we decompose:

$$\chi_{Q^{\tilde{m}}_{j}}S_{j, \leq 2^{-i}} \left( \sum_{\xi} \sum_{k} g_{\xi, \xi^{2}+k} \right) \left( \sum_{\eta \in \Theta_{p}} f_{\eta, \leq 2^{-i}} \right)=$$

$$\sum_{\a} \sum_{l} \sum_{m \in B^{\a}} \sum_{m' \in A^{\a}} \sum_{k} \sum_{\xi} \sum_{\eta \in \Theta_{p}} \chi_{Q_{i}^{m}} g_{\xi, \xi^{2}+k} \cdot \chi_{T_{\eta}^{m',l}} \chi_{Q^{\tilde{m}}_{j}} f_{\eta, \leq 2^{-i}}$$  

Then we continue the exact same argument, just that we always replace $f_{\eta, \leq 2^{-i}}$ by $\chi_{Q^{\tilde{m}}_{j}} f_{\eta, \leq 2^{-i}}$. We end up with:

$$\sum_{\eta} || \chi_{Q^{\tilde{m}}_{j}} (g_{i, \leq 2^{i-2}} \cdot f_{j, \leq 2^{-i}})_{\eta, \leq 2^{-i}}||^{2}_{\mathcal{Y}_{\eta}} \les 2^{i-j}  ||g_{i,\leq 2^{i-2}}||_{\mathcal{D}X^{0,\q,1}} \left( \sum_{\eta} ||\chi_{Q^{\tilde{m}}_{j}} f_{\eta, \leq 2^{-i}}||^{2}_{Y_{\eta}} \right)^{\q} $$

Summing over the set of $\tilde{m}$ with the property $Q_{\tilde{m}}^{j} \cap L^{k}_{i} \ne \emptyset$, and then taking the suppremum with respect to $k$ and $L$ gives us (\ref{a5}).

\end{proof}

\vspace{.1in}

\subsection{Estimates: {\mathversion{bold}$ X^{0,\q,1}_{i} \cdot Y_{j, \leq 2^{-i}} \rightarrow X^{0,\q,1}_{j, \geq 2^{-i}}$}}

\noindent
\vspace{.1in}

The main estimate in this section is the following:

\begin{p3} We have the estimate

\beq \label{a17}
|| v_{j, \leq 2^{-i}} \cdot  u_{i}||_{X^{0,-\q}_{j, \geq 2^{-i}}} \les 2^{i-j} i^{\q} ||v_{j, \leq 2^{-i}} ||_{Y_{j}} \cdot ||u_{i}||_{X^{0,\q,1}}
\eeq

\end{p3}

The proof of this result is split again into two parts. We claim:

\beq \label{a9}
|| v_{j, \leq 2^{-i}} \cdot  u_{i, \leq 2^{i-2}}||_{X^{0,-\q}_{j, \geq 2^{-i}}} \les 2^{i-j} i^{\q} ||v_{j, \leq 2^{-i}} ||_{Y_{j}} \cdot ||u_{i, \leq 2^{i-2}}||_{X^{0,\q,1}}
\eeq

\beq \label{a10}
|| v_{j, \leq 2^{-i}} \cdot  u_{i, \geq 2^{i-2}}||_{X^{0,-\q}_{j, \geq 2^{-i}}} \les 2^{i-j} i^{\q} ||v_{j, \leq 2^{-i}} ||_{Y_{j}} \cdot ||u_{i, \geq 2^{i-2}}||_{L^{2}}
\eeq

\begin{proof} We decompose $v_{j, \leq 2^{-i}}$ as in (\ref{dec1}) and  $u_{i, \leq 2^{i-2}}$ as in (\ref{dec2}). We know from previous section that:

\beq \label{k7}
||g_{i, \leq 2^{i-2}} \cdot f_{j, \leq 2^{-i}}||^{2}_{X^{0,\q,1}_{j, \geq 2^{-i}}} \approx \sum_{p} ||g_{i, \leq 2^{i-2}} \cdot \sum_{\eta \in \Theta_{p}} f_{\eta, \leq 2^{-i}}||^{2}_{X^{0,\q,1}_{j, \geq 2^{-i}}}
\eeq

For a fixed $d_{2} \geq 2^{-i}$, we decompose:

$$S_{j, d_{2}} (g_{i,\leq 2^{-i}} \cdot \sum_{\eta \in \Theta_{p}} f_{\eta, \leq 2^{-i}})= S_{j, d_{2}} \left( \sum_{k=0}^{2^{2i-2}} \sum_{\xi} g_{\xi, \xi^{2}+k} \right) \left( \sum_{\eta \in \Theta_{p}} f_{\eta, \leq 2^{-i}} \right) $$

From Lemma \ref{g1} we know that $g_{\xi, \xi^{2}+k} \cdot f_{\eta, \leq 2^{-i}}$ is supported in $A_{j, d_{2}}$ iff $|\cos{\a}| \leq 2^{-i} d_{2}$, where $\a$ is the angle between $\xi$ and $\eta$. The angle between any two $\eta$'s in $\Theta_{p}$ is at most $2^{i-j} \leq 2^{-2i}$ and the angle between any two $\xi$'s in $\Xi^{i}$ is either at least $2^{-2i}$ or the same. Therefore the $\xi's$ involved in the above summation have an angular localization in a set of cardinality $\approx 2^{i}d_{2}$; we just keep this in mind and not formalize it. What is really important is that we sum over a set containing $\approx 2^{3i}d_{2}$ $\xi$'s. 

For each $g_{\xi,\xi^{2}+k}$ and $f_{\eta, \leq 2^{-i}}$ we can apply the result in (\ref{m9}):

$$|| f_{\eta, \leq 2^{-i}} \cdot  g_{\xi,\xi^{2}+k}||_{L^{2}} \les  2^{-\frac{i+j}{2}} ||f_{\eta, \leq 2^{-i}}||_{Y_{j}} ||g_{\xi,\xi^{2}+k}||_{L^{2}}$$

Using the result in Lemma \ref{g1} we can perform the summation with respect to $\eta \in \Theta_{p}$:

 $$||\sum_{\eta \in \Theta_{p}} f_{\eta, \leq 2^{-i}} \cdot  g_{\xi,\xi^{2}+k}||_{L^{2}} \les  2^{-\frac{i+j}{2}} ||\sum_{\eta \in \Theta_{p}} f_{\eta, \leq 2^{-i}}||_{Y_{j}} ||g_{\xi,\xi^{2}+k}||_{L^{2}}$$

Then we can perform the summation with respect to $\xi$:

$$||\sum_{\eta \in \Theta_{p}} f_{\eta, \leq 2^{-i}} \cdot  g_{i, 2^{-i}k}||_{L^{2}_{j,d_{2}}} \les  2^{-\frac{i+j}{2}} (2^{3i}d_{2})^{\q} ||\sum_{\eta \in \Theta_{p}} f_{\eta, \leq 2^{-i}}||_{Y_{j}} ||g_{i, 2^{-i}k}||_{L^{2}}$$

\noindent
followed by the one with respect to $k$:

$$||\sum_{\eta \in \Theta_{p}} f_{\eta, \leq 2^{-i}} \cdot  g_{i}||_{L^{2}_{j,d_{2}}} \les  2^{-\frac{i+j}{2}} (2^{3i}d_{2})^{\q} ||\sum_{\eta \in \Theta_{p}} f_{\eta, \leq 2^{-i}}||_{Y_{j}} ||g_{i}||_{X^{0,\q,1}}$$

In the end we perform the summation with respect to $p$ and pass to $X^{0,\q}$ norm:

$$|| f_{j, \leq 2^{-i}} \cdot  g_{i}||_{X^{0,\q}_{j,d_{2}}} \les  2^{i-j}  || f_{j, \leq 2^{-i}}||_{Y_{j}} ||g_{i}||_{X^{0,\q,1}}$$

We sum up with respect to $d_{2}$ (over a set of cardinality $\approx 2i$) to obtain the statement in (\ref{a9}).

Now we continue with the proof of (\ref{a10}). The approach is similar to the one above, but we still outline the main steps. We decompose $v_{j, \leq 2^{-i}}$ as in (\ref{dec1}) and  $u_{i, \geq 2^{i-2}}$ as in (\ref{dec2}). We know from previous section that:

\beq \label{k8}
||g_{i, \geq 2^{i-2}} \cdot f_{j, \leq 2^{-i}}||^{2}_{X^{0,\q,1}_{j, \geq 2^{-i}}} \approx \sum_{p} ||g_{i, \geq 2^{i-2}} \cdot \sum_{\eta \in \Theta_{p}} f_{\eta, \leq 2^{-i}}||^{2}_{X^{0,\q,1}_{j, \geq 2^{-i}}}
\eeq

For a fixed $d_{2} \geq 2^{-i}$, we decompose:

$$S_{j, d_{2}} (g_{i,\leq 2^{-i}} \cdot \sum_{\eta \in \Theta_{p}} f_{\eta, \leq 2^{-i}})= S_{j, d_{2}} \left( \sum_{l} \sum_{n} \sum_{\xi \in \Xi^{i}_{n}} g_{\xi, l} \right) \left( \sum_{\eta \in \Theta_{p}} f_{\eta, \leq 2^{-i}} \right) $$

From Lemma \ref{g1} we know that $g_{\xi, l} \cdot f_{\eta, \leq 2^{-i}}$ is supported in $A_{j, d_{2}}$ iff $|\cos{\a}| \leq |\xi|^{-1} d_{2}$, where $\a$ is the angle between $\xi$ and $\eta$. The angle between any two $\eta$'s in $\Theta_{p}$ is at most $2^{i-j} \leq 2^{-2i}$ and the angle between any two $\xi$'s in $\Xi^{i}_{n}$ is either at least $n^{-1}=|\xi|^{-1}2^{-i}$ or the same. Therefore the $\xi's$ involved in the above summation have an angular localization in a set of cardinality $\approx 2^{i}d_{2}$; we will just keep this in mind and not formalize it. What will be really important is that we sum over a set containing $\approx 2^{i}d_{2}$ $\xi$'s from $\Xi^{i}_{n}$ and over a set containing $\approx 2^{3i}d_{2}$ $\xi$'s from $\Xi^{i}$. 

The setup is exactly like in the proof of (\ref{a9}) with $l$ playing the role of $k$ and the proof can be continued in the same fashion. 

\end{proof}

\vspace{.1in}

\subsection{Estimates: {\mathversion{bold}$ X^{0,\q,1}_{i} \cdot X_{j, \geq 2^{-i}}^{0,\q,1} \rightarrow \mathcal{Y}_{j,\leq 2^{-i}}$}}

\noindent
\vspace{.1in}

The main estimate in this section is the following:

\begin{p3} We have the estimate

\beq \label{a19}
|| v_{j, \geq 2^{-i}} \cdot  u_{i}||_{\mathcal{Y}_{j, \leq 2^{-i}}} \les 2^{i-j} i^{\q} ||v_{j, \geq 2^{-i}} ||_{X^{0,\q,1}_{j,\geq 2^{-i}}} \cdot ||u_{i}||_{X^{0,\q,1}}
\eeq

\end{p3}

This can be obtained by duality from {\mathversion{bold}$ X^{0,\q,1}_{i} \cdot Y_{j, \leq 2^{-i}} \rightarrow X_{j,\geq 2^{-i}}^{0,\q,1}$}.

\subsection{Bilinear estimates on dyadic regions}

\begin{proof}[Proof of Theorem \ref{tb2}] We decompose

$$B(u_{i},v_{j}) = S_{j, \leq 2^{-i}}B(u_{i},v_{j}) +  S_{j, \geq 2^{-i}}B(u_{i},v_{j}) =$$

$$ S_{j, \leq 2^{-i}}B(u_{i},v_{j, \leq 2^{-i}}) + S_{j, \leq 2^{-i}}B(u_{i},v_{j, \geq 2^{-i}}) + $$

$$ S_{j, \geq 2^{-i}}B(u_{i},v_{j, \leq 2^{-i}}) + S_{j, \geq 2^{-i}}B(u_{i},v_{j, \geq 2^{-i}})$$

For the first term we make use of (\ref{a3}) to obtain:

$$|| S_{j, \leq 2^{-i}}B(u_{i},v_{j, \leq 2^{-i}})||_{\mathcal{Y}^{s}} \approx 2^{sj} || S_{j, \leq 2^{-i}}B(u_{i},v_{j, \leq 2^{-i}})||_{\mathcal{Y}} \les $$

$$i^{\q} 2^{i} 2^{(s-1)j} ||\nabla v_{j, \leq 2^{-i}}||_{Y} ||\nabla u_{i}||_{\mathcal{D}X^{0,\q}} \les $$

$$i^{\q} 2^{2i} 2^{sj} ||v_{j, \leq 2^{-i}}||_{Y} ||u_{i}||_{\mathcal{D}X^{0,\q}} \approx $$

$$i^{\q} 2^{(2-s)i} ||v_{j, \leq 2^{-i}}||_{Y^{s}} ||u_{i}||_{\mathcal{D}X^{s,\q}} \les i 2^{(2-s)i} ||v_{j, \leq 2^{-i}}||_{Y^{s}} ||u_{i}||_{\mathcal{D}Z^{s}}$$

For the second term we make use of (\ref{a19}) to obtain:

$$||S_{j, \leq 2^{-i}}B(u_{i},v_{j, \geq 2^{-i}})||_{\mathcal{Y}^{s}} \approx 2^{sj} ||S_{j, \geq 2^{-i}}B(u_{i},v_{j, \geq 2^{-i}})||_{\mathcal{Y}} \les $$

$$i^{\q} 2^{i} 2^{(s-1)j} ||\nabla v_{j, \geq 2^{-i}}||_{X^{0,\q}} \cdot ||\nabla u_{i}||_{\mathcal{D}X^{0,\q}} \les $$

$$i^{\q}  2^{(2-s)i} ||v_{j, \geq 2^{-i}}||_{X^{s,\q}} \cdot || u_{i}||_{\mathcal{D}X^{s,\q}} \les $$

$$i^{\frac{3}{2}}  2^{(2-s)i} ||v_{j, \geq 2^{-i}}||_{Z^{s}} \cdot || u_{i}||_{\mathcal{D}Z^{s}} $$

For the third term we use of (\ref{a17}) to obtain:

$$||S_{j, \geq 2^{-i}}B(u_{i},v_{j, \leq 2^{-i}})||_{X^{s,-\q,1}} \approx 2^{sj} ||S_{j, \geq 2^{-i}}B(u_{i},v_{j, \leq 2^{-i}})||_{X^{0,-\q,1}} \les $$

$$i^{\frac{3}{2}}2^{i} 2^{(s-1)j} ||\nabla v_{j, \leq 2^{-i}}||_{Y} ||\nabla u_{i}||_{\mathcal{D}X^{0,\q}} \les $$

$$i^{\frac{3}{2}} 2^{(2-s)i} ||v_{j, \leq 2^{-i}}||_{Y^{s}} ||u_{i}||_{\mathcal{D}X^{s,\q}} \les i^{2} 2^{(2-s)i} ||v_{j, \leq 2^{-i}}||_{Y^{s}} ||u_{i}||_{\mathcal{D}Z^{s}} $$

The fourth term had been handled in Theorem \ref{bil}. By adding all the estimates we obtain:

$$||B(u,v)||_{W^{s}_{j}} \les i^{\frac{3}{2}} 2^{(2-s)i} ||u||_{\mathcal{D}Z^{s}_{i}} ||v||_{Z^{s}_{j}}$$

In the end we can recover the decay via an argument similar to the one in Proposition \ref{pdec}, part a). One would notice that over there we had to recover decay of type $\mathcal{D}_{j}$ and all we used is that the high frequency comes with that decay. We already worked out the conservation of decay for the first term, see (\ref{a4}).

\end{proof}

\section{Bilinear estimates - Proof of Theorem \ref{bg}}

This is a standard argument once we have the bilinear estimates on dyadic pieces, see (\ref{be11}), (\ref{be22}) and (\ref{b7}). For reference one could use Part 1, see the corresponding section there.

\end{document}